\documentclass[12pt]{amsart}
\usepackage{amsmath,amssymb,enumerate, graphicx}
\usepackage[frenchb]{babel}
\usepackage[latin1]{inputenc}
\providecommand{\Phi}{\textbf{\Phi}}

\renewcommand{\L}{\mathcal{L}}

\renewcommand{\P}{\mathbb{P}}

\providecommand{\sup}{\mathrm{sup}}
\providecommand{\inf}{\mathrm{inf}}
\providecommand{\dim}{\mathrm{dim} \,}
\providecommand{\C}{\mathcal{C}} 
\providecommand{\L}{\mathcal{L}}

\providecommand{\G}{\mathcal{G}} \providecommand{\A}{\mathcal{A}}
\providecommand{\B}{\mathcal{B} } 
\providecommand{\R}{\mathbb{R}}
 \providecommand{\E}{\mathbb{E}}
\providecommand{\P}{\mathbb{P}} \providecommand{\N}{\mathbb{N}}
\providecommand{\Z}{\mathbb{Z}}

\newtheorem{theo}{Théorème}[section]

\newtheorem{lem}[theo]{Lemme}
\newtheorem{souslem}[theo]{Sous lemme}
\newtheorem{prop}[theo]{Proposition}
\newtheorem{rem}[theo]{Remarque}

\newtheorem{defi}[theo]{Définition}
\begin{document}

\title[Nombre de points visit\'es sur un amas de percolation  ]{\textbf{
Sur le nombre de points visit\'es  par une marche al\'eatoire sur 
un amas infini de percolation.}}
\author[Cl\'ement Rau]{Clément
Rau\protect\footnotemark{$^1$} }
\date{}

\begin{abstract}

Dans cet article, on s'int\'eresse à une marche al\'eatoire simple
sur un amas infini issu d'un processus de  percolation  surcritique sur les ar\^etes de 
 $\Z^d \  (d \geq 2)$ de loi $Q$. On montre que la transform\'ee de Laplace du nombre 
 de points visit\'es au temps $n$, not\'e $N_n$, a un comportement similaire
 au cas o\`u la marche \'evolue dans $\Z^d$. Plus précisément,
 on établit que  pour tout $0<\alpha<1$, il existe des constantes 
 $C_i, \ C_s >0$ telles que  pour presque toute réalisation de la
 percolation telle que l'origine appartienne à l'amas infini  et  pour $n$ assez grand,
 $$  e^{-C_i n^{ \frac{d}{d+2} } }  \leq \E_0^{\omega} ( \alpha^{N_n} )
 \leq e^{-C_sn^{ \frac{d}{d+2} }}.$$
Le point principal du travail réside dans l'obtention de la borne supérieure.
Notre approche consiste  dans un premier temps, à  trouver  une inégalité 
isopérimétrique sur l'amas
infini, et dans un deuxi-\\
-ème temps à la remonter  sur un produit en couronne,
 ce qui nous permet 
alors d'obtenir une majoration de la probabilité de retour d'une certaine marche
sur ce produit en couronne. L'introduction d'un produit en couronne est 
justement motivée par le fait que la probabilité de retour sur un tel graphe 
s'interprète comme l'espérance 
de la transformée de Laplace du nombre de points 
visités. 

\end{abstract}

\maketitle 

\section{Introduction et résultats}

\protect\footnotetext[1]{ CMI, 39 rue Joliot-Curie,
13013 Marseille, rau@cmi.univ-mrs.fr  }

Soit $d\geq 2$.
On appelle percolation de Bernoulli de paramètre $p$ le
sous graphe aléatoire de la grille de dimension $d$
obtenu en supprimant (resp. gardant) une ar\^ete avec
probabilité $p$ (resp. $1-p$) de fa\c{c}on indépendante
pour les différentes ar\^etes. On notera $\omega$ une réalisation typique de la percolation.
On appelle alors amas infini une composante connexe
infinie du graphe $\omega$.
On montre qu'une telle composante connexe infinie existe
et est presque s\^urement unique
si le paramètre $p$ est choisi au dessus d'une certaine valeur critique $p_c$.
Une construction plus formelle de la percolation est
donnée plus loin dans cette introduction.

La percolation est un modéle important de la mécanique
statistique des milieux desordonnés,
un 'concept unificateur' pour reprendre l'expression de
P.G. De Gennes dans \cite{gennes},
qui intervient aussi dans de nombreuses applications, par
exemple dans les problèmes de diffusion
dans un environnement non homogène que l'on rencontre
dans la recherche pétrolière.
Depuis son introduction en 1956 par J.M. Hammersley,  la
percolation a également donné lieu
a une jolie théorie mathématique qui recèle encore bien
des défis.
Nous renvoyons aux livres de Kesten et Grimmett (voir
\cite{kesten} et \cite{grimmett}) pour une introduction
aux outils mathématiques de la percolation. On y trouvera
en particulier de nombreux résultats
sur la géométrie des amas infinis.

Depuis quelques années, différents auteurs se sont
attachés à dévelop-per
la théorie du potentiel des amas infinis ou, en d'autres
termes, à décrire le comportement d'une
marche aléatoire évoluant sur un amas infini
(La 'fourmi dans un labyrinthe' pour reprendre une autre
expression de P.G. De Gennes).
Les premières bornes sur le noyau de la chaleur sur un
amas infini ont été demontrées par P. Mathieu et E. Rémy
(voir \cite{pierre}) à l'aide d'estimées
du profil isopérimétrique d'un amas infini. Puis, M. Barlow
(voir \cite{barjo}) a obtenu des estimées de type gaussien.
Ces premiers résultats ont ensuite permis de prouver la
convergence de la marche aléatoire vers un mouvement
brownien sous la forme d'un principe d'invariance valable
pour presque toute réalisation de $\omega$, voir \cite{SS}
\cite{BB} et \cite{MP}.
L'objet du travail présenté ici est de compléter ce
panorama en donnant des estimées précises
sur le nombre de points visités par la marche aléatoire
simple symétrique évoluant sur un amas infini.

Le processus de percolation est défini de la manière
suivante. Pour $d\geq 2$, notons $E^d$ l'ensemble des ar\^etes de $\Z^d$ défini par,
$$E^d = \{ (x,y) ;
\sum_{i=1..d} |x_i -y_i| =1 \},$$
o\`u $x=(x_1,...,x_d)$ et $y=(y_1,...,y_d)$. Pour $p\in ]0,1]$, soit  $\omega$
  le sous graphe aléatoire de  $\L^d:=(\Z^d,E^d)$
 obtenu en gardant (resp. effa\c{c}ant) une ar\^ete   avec  probabilité $p$
  (resp. $1-p$)  de manière indépendante pour les différentes ar\^etes de $E^d$.
  On identifie  ce  sous graphe de $\Z^d$  avec l'application 
  $ \omega :E^d \rightarrow \{ 0,1\}$
telle que  $\omega (x,y)=1  $ si l'ar\^ete  $(x,y) $ est présente dans 
 $\omega$ (on dira qu'une telle ar\^ete est 'ouverte') et  
 $\omega(x,y) =0 $ sinon. On munit  $\{0,1\}^{E^d}$ de la mesure de probabilité
 $Q$   sous laquelle les variables aléatoires
  $(\omega (e ), e\in \E^d)$ sont indépendantes et suivent des lois de  
  Bernoulli(p). 
  Soit $\C $ 
 la composante connexe de $\omega$ contenant l'origine, $|\C|$ son cardinal et $p_c$ la probabilité 
 critique,
  $$p_c =\sup \{ p ; Q ( |\C |= +
\infty) =0
\}.$$ 
On sait que $0<p_c<1$(voir \cite{grimmett}) , et désormais, on suppose que  $p>p_c$. On se place sur 
l'événement 
$\{|\C|=+\infty\}$, $\C$ est alors l'unique amas infini. On
considère alors sur $\C$ la marche aléatoire 
suivante: $X_0=x $  et 
$X_{n+1}  $ est choisi uniformément parmi les voisins de 
$X_n$ dans $\C$. 
$$\text{ie: }\ 
\P(X_{n+1}=y|X_n=x)=\frac{\omega(x,y)}{\underset{z; 
(x,z)\in E^d}\sum
\omega(x,z)}.$$
 Le but de cet article  est d'estimer le nombre de points visités par la marche $X$. Plus 
précisément, posons 
 $N_n= |\{X_0, X_1,...,X_n \}|$,  $\P_x^{\omega} $ la loi de la marche
 issue de $x$,  et 
   $\E_x^{\omega} $ son espérance.
   Le principal résultat est:
  \begin{theo} \label{theoinitial} Pour tout  $\alpha \in ]0,1[$, il existe 
  deux constantes 
  $C_i, C_s >0$ telles que  
 Q p.s  sur l'événement  $|\C|=+\infty$, et pour  $n$ assez grand: 
 $$  e^{-C_i n^{\frac{d}{d+2}}}   \leq \mathbb{E}_0^{\omega}(\alpha^{N_{n}} ) 
 \leq e^{-C_s n^{\frac{d}{d+2}}}.$$
 \end{theo}
 \begin{rem} Ce résultat est également valide
 pour la marche aléatoire à temps continu qui
 attend un temps exponentiel entre chaque
 saut.
 \end{rem}
 Dans le cas o\`u $\omega=\L^d$ cette expression a déjà été
 étudiée par 
   Donsker M D et Varadhan S R S (voir  \cite{ donsk}). Ils 
 prouvent, en particulier pour la marche aléatoire  simple sur la grille
 $\Z^d$, le 
 théorème suivant:
 \begin{theo} (voir \cite{ donsk})
Il existe une constante $ c(d,\alpha) >0$ telle que, $$   \underset{n\rightarrow \infty }{lim} \frac{1}
{n^{\frac{d}{d+2} }}\log\E_0^{\Z^d}(\alpha^{N_n} ) =-c(d,\alpha)
$$
\end{theo} 
Ce résultat est prouvé par double inégalité. La  strat\'egie adoptée dans la
preuve de Donsker M D et Varadhan S R S, pour obtenir une majoration  de 
l'expression $\underset{n\rightarrow \infty }{lim} \frac{1}
{n^{\frac{d}{d+2} }}\log\E_0^{\Z^d}(\alpha^{N_n} )$, ne semble pas 
se généraliser dans un amas infini de percolation pour
plusieurs raisons. 
Par exemple, la symétrie de  $\Z^d$  est un point crucial 
dans leur preuve, qui n'est évidemment pas satisfait dans un amas.
En particulier,  la marche aléatoire sur $\Z^d$ satisfait
une propriété de martingale qui n'est plus vraie pour la
marche aléatoire sur l'amas de percolation. 

La  méthode  developpée ici  repose sur le fait suivant: 
$\E_0^{\omega}(\alpha^{N_{n}} )$ peut s'écrire 
comme une probabilité de retour à l'origine d'une marche $Z$ construite 
à partir de $X$ dans un graphe plus "gros" que  $\C$, qui sera un produit 
en couronne.  Trouver une borne supérieure de $\E_0^{\omega}(\alpha^{N_{n}} )$ 
revient donc à trouver une borne supérieure de la probabilité de retour de la
marche $Z$. On sait que  les
inégalités isopérimétriques sont un outil important pour prouver des
inégalités fonctionnelles comme celles de Poincaré ou de  Nash,  qui elles m\^emes,
permettent d'obtenir des bornes supérieures du noyaux d'une marche simple
(voir \cite{coulhon}). 
 On \'etudie donc le profil isopérimétrique
 sur ce produit en couronne. Gr\^ace aux  récents travaux
 d'A.Erschler,
 on sait contr\^oler l'isopérimétrie du produit en couronne de deux graphes à
 partir de l'isopérimétrie de chacun d'entre eux.
Ici, un des deux graphes étant le graphe ayant comme ensemble de points $\C$,
on est finalement ramené à étudier de manière assez fine 
la géométrie d'un amas et  ses propriétés isopérimétriques.  
Notons $\B_n =[-n;n]^d$ et  $\C_n$ la composante
connexe  de  $\C\cap \B_n$ contenant l'origine.  Avec des techniques similaires 
à celles de  \cite{pierre}, 
on prouve la propriété suivante.
\begin{prop} \label{is de base} Soient  $\gamma>0$ et   $p>p_c$. 
 Il existe  $\beta>0 $  tel que pour tout $c>0$ :
 $Q$ p.s  sur  $\#\C=+\infty$, pour  $n$ assez grand, on ait,
 \begin{eqnarray} \label{betac}
  \frac{|\partial_{\C^g} A | }{f_c(|A|) } \geq
\beta \ \ \ \text{pour tout sous-ensemble $A$ connexe de } \C_n, 
\end{eqnarray}
o\`u
 $ f_c(x)=
 \left\lbrace
\begin{array}{l}
  1  \ \ \ \ \ \  \  si \  x <c  n^{\gamma} \\
 x^{1-\frac{1}{d}}    \ \  si\  x \geq c n^{\gamma},\\
\end{array}
\right.$ \\
et  $\partial_{\C^g} A =
\{  (x,y)\in \E^d ; \ \omega(x,y)=1 \ \text{et} \ x\in A \ ; \ y \not\in A \}.$\\
Le rang à partir duquel l'inégalité (\ref{betac}) est satisfaite dépend de l'amas
$\omega$ et de $c$. 
\end{prop}

\begin{rem} Les m\^emes techniques
s'appliquent  pour des marches aléatoires sur $\Z^d$ aux plus
proches voisins,  symétriques et dont les taux de transitions sont bornés
supérieurement et inférieurement (condition d'ellipticité).
Cette dernière condition n'est pas satisfaite pour la
percolation.
 \end{rem}

Après avoir fixé quelques notations, ce papier se décompose en  4 parties. 
Dans la partie 2, nous définissons  un produit en couronne, nous expliquons
l'int\^eret d'un tel outil dans notre cas et nous traitons le problème de 
l'isopérimétrie sur un produit en
couronne. Dans la partie 3, nous étudions le profil isopérimétrique d'un amas et
nous prouvons en particulier la propositions \ref{is de base}. A l'aide des
parties 2 et 3, nous établissons dans la partie 4, une borne supérieure pour la probabilité de
retour à l'origine dans le produit en couronne, et nous en déduisons  une borne supérieure
pour $\E_0^{\omega} (\alpha^{N_n})$.
Enfin, dans la partie 5 nous prouvons la borne inférieure du théorème
\ref{theoinitial}.

\subsection{Notations}$\ $ \\
\begin{enumerate} [$\bullet$]
\item  Nous utiliserons le symbole  $:= $  pour définir une nouvelle 
quantité et on notera  $|A|$ ou bien $\# A$ le cardinal d'un ensemble $A$.
\item La somme  disjointe de deux ensembles $A$ et $B$ quelconques sera notée 
$A \dot{\bigcup} B $. Cette notation représente un ensemble  en bijection avec
$(A\times \{0\} )\cup( B\times \{1\}).$   
\item  On notera  $ D(0,x)\ $ la longueur minimale d'un
chemin constitué d'ar\^et-es ouvertes reliant 
$x$ à $0$,
(appelé aussi distance chimique) et  on pose  
 $B_r(\C)=\{x\in\Z^d;\ D(0,x)\leq r\}$.  
On utilisera également la norme  $N_1$ sur  $\Z^d$ définie par 
$N_1(x)=\sum_{i=1...d}
 |x_i|,$ si $x=(x_1,...,x_d).$
  \item $C$ et $  c$  représenteront des constantes dont la valeur
  peut évoluer de ligne en ligne mais la dépendance en $n$ ou  $\omega$
  sera spécifiée par l'ordre des quantificateurs. La constante $\beta$ 
  dépendra uniquement de  $p$.
\item Un graphe $G$ est un couple $(V(G),E(G))$, o\`u  $V(G)$ désigne l'ensemble
des points  de $G$ (vertices of $G$) et $E(G)$ désigne l'ensemble des ar\^etes 
de $G$. Dans ce papier, les graphes seront non orientés (
sauf mention explicite
comme par exemple dans le lemme \ref{soushypergraphe}), et o\`u les ar\^etes 
de type $(x,x)$ sont exclues.
Un sous graphe de $G$  est un graphe  $G'$ tel que  $V(G')\subset V(G)$ et 
 $E(G') \subset E(G)$. 
 \item On utilisera en particulier le graphe $\C^g$ défini par 
 $V(\C^g)=\C$ et $E(\C^g)=\{(x,y)\in E^d;\ x,y\in\C \text{ et } \omega(x,y)=1\}$ et le graphe 
  $\C_n^g$ défini par $V(\C_n^g)=\C_n$ (on rappelle que $C_n$ désigne la
  composante connexe de $\C\cap[-n,n]^d$ contenant l'origine) et
   $E(\C_n^g)=\{(x,y)\in E^d;\ x,y\in\C_n \ \text{et}\ \omega(x,y)=1\}$.
   De fa\c{c}on générale, lorsqu'il n y a pas d'ambiguité
   sur les ar\^etes,  $V^g$ désignera le graphe  ayant cet
   ensemble d'ar\^etes et dont    l'ensemble des points  est
   $V$ . On notera \'egalement $\B_n^g$ le graphe ayant $\B_n=[-n,n]^d$ comme ensemble
   de points et o\`u $E(\B_n^g)=\{ (x,y)\in E^d;\ x \ \text{et} \ y \in \B_n \}.$
\item  Soit $G$ un graphe, pour $A\subset V(G)$, 
on notera
  $$ \partial_G A=\{ (x,y)\in E(G);  \ x \in A  \ \text{et} \ y\in V(G)-A    \}.$$
  Cette notation est evidemment cohérente avec la définition du bord 
  $\partial_{\C^g}$ dans la
  proposition \ref{is de base}. 
En notant $\L^d=(\Z^d,E^d)$,  le bord "classique" dans  $\Z^d$ d'un
ensemble $A$ est donc noté $\partial_{\L^d} A$.
\item On travaillera plut\^ot avec les fonctions de F\o lner lorsqu'on parlera
d'isopérimétrie. Soit $G$ un graphe,  on d\'esigne par $Fol_G$ la fonction
de F\o lner de $G$, d\'efinie par:
$$ Fol_G(k) = \min \{ |U|; \ U\subset V(G)  \ \text{et} \  
 \frac{|\partial_{G} U|}{ |U|}\leq \frac{1}{k} \}.$$
$\bullet$ Si $G'$ est un sous graphe de $G$,  on notera:
$$ Fol_{G'}^{G}(k) = \min \{ |U|; \ U\subset V(G') \ \text{et} \  
\frac{|\partial_{G} U|}{ |U|}\leq \frac{1}{k} \}.$$
$Fol_{G'}^{G}$ est toujours une fonction de F\o lner de $G'$ mais o\`u le bord est
compt\'e dans $G$.

\end{enumerate}
\section{Produit en couronne.}
Dans cette section, après avoir donné la définition d'un produit en couronne,
 on motive l'introduction  de tels graphes par la propriété
  \ref{proba-Nn}, en consid\'erant une marche aléatoire construite à partir
  des noyaux de transition de la marche $X$, sur un certain produit en couronne.
   Enfin, on estime la fonction de F\o lner du produit
en couronne en question à l'aide des fonctions de F\o lner des graphes considérés.
\subsection{Définition } \ 
Soient $A$ et $B$ deux graphes, et $b_0$ un point de $V(B)$. Pour une fonction
$f:V(A)\rightarrow V(B)$, on appelle support de $f$, l'ensemble 
$\{a\in V(A);\ f(a)\neq b_0 \}.$ 
\begin{defi} Le produit en couronne 
 $A \wr B$  de deux graphes $A$ et $B$, est le graphe suivant:  \\
- $V(A\wr B)$ est l'ensemble des couples   $(a,f)$ o\`u  $a\in V(A)$  et   $f: V(A)
\rightarrow V(B)$  est à support fini, \\
-les ar\^etes sont définies de la manière suivante:  
$ \Bigl((a_1,f_1) ,(a_2,f_2)\Bigr) \in E(A\wr B))$ si et seulement si
\begin{eqnarray*}
  \begin{cases}
  \;  a_1 =a_2 \  \text{et} \ \forall x\neq a_1 \  f_1(x)=f_2(x) \  \text{et} \
  ( f_1(a_1), f_2(a_1))\in E(B) , \\
\;\text{ou}\\
\;f_1=f_2  \ \ \text{et} \ \  (a_1,  a_2)\in E(A).
\end{cases}
\end{eqnarray*}
\end{defi}
 Si $(a,f) \in V(A\wr B)$, l'élément $f$ est appelée la configuration.
 On appellera le graphe
 $A$  la "base" du produit en couronne $A\wr B$.
\begin{rem}  Si $A$ et $B $ sont des graphes de  Cayley  de groupes, le 
produit en couronne de  $A$ et  $B$ est le graphe de  Cayley  du produit 
en couronne de ces groupes, avec l'ensemble "standard" de générateurs 
construits à partir des générateurs de $A$ et $B$ (voir \cite{ersh}). [
Le produit en couronne de deux groupes $A$ et $B$ est le produit semi direct 
de  $A$  et $\sum_A B$ o\`u  $A$ agit sur  $\sum_A B $ par 
 $\  ^{a}f(x)=f(  xa^{-1} )$. ] 
\end{rem}

\subsection{ Marches aléatoires}
\label{RW}
Dans notre cas, on prend  $A=\C^g$ et   pour $B$  le graphe de  Cayley  du
groupe  $(\frac{\Z}{2\Z},+) $ avec  $\overline{1}$  comme générateur. 
On notera 'encore'  $\frac{\Z}{2\Z}$ ce graphe. On choisit
 $b_0=\overline{0}$. Soit $o$ le point de $ A\wr B$ tel que $o=(0,f_0)$
 o\`u $0$  est l'origine de $\Z^d$ et  $f_0$  est la configuration 
 qui vaut  $\overline{0}$ en tous points. Notons $p(\ ,\ )$ les noyaux de
 transition de la marche  $X$  définie dans l'introduction, ie: pour tout 
 $a,b \in \C \ \  p(a,b)=\P_x^{\omega}(X_1=y)=\frac{\omega(a,b)}{\nu(a)}$ o\`u $\nu(a)$ est le nombre
 de voisins dans $\C$ de $a$. Pour $\alpha \in ]0;1[ $, on consid\`ere  alors la
  marche  $Z$ sur  $\C^g \wr \frac{\Z}{2\Z}$ définie par  $Z_0=o $ et 
dont les noyaux de   probabilité sont:
\begin{eqnarray*}
\tilde{p} \Bigl( (a,f),(b,f_{\substack{a,\bar{0} \\ b,\bar{0} } } )\Bigr) = \alpha^2 p(a,b), \\ 
 \tilde{p} \Bigl((a,f),(b,f_{\substack{a,\bar{1} \\ b,\bar{0}  } } )\Bigr) =
\alpha(1-\alpha) p(a,b), 
 \\ \tilde{p} \Bigl( (a,f),(b,f_{\substack{a,\bar{0} \\b,\bar{1} }} )\Bigr) = \alpha(1-\alpha) p(a,b),\\ 
 \tilde{p} \Bigl(
(a,f),(b,f_{\substack{a,\bar{1} \\ b,\bar{1}}} )\Bigr) = (1-\alpha)^2 p(a,b),
\end{eqnarray*}
o\`u \begin{eqnarray*} f_{ \substack{a,x \\ b,y} }: &u& \mapsto  f(u)  
\text{ pour $ u \not=  a $ ou $ b,$}\\
 &a&\mapsto x, \\
&b& \mapsto y. 
\end{eqnarray*}
L'interprétation  de cette marche est la
suivante:  imaginons qu'il y ait une
lampe en chaque point de $\C$, qui soit allumée [resp. éteinte] lorsque la configuration en ce
point vaut $\overline{1}$  [resp. $\overline{0}$]. Supposons maintenant qu'à
un certain instant le marcheur se trouve en un certain point de $\C$. En une 
unité de temps, il allume [resp. éteint] la lampe où il se trouve avec probabilité 
 $1-\alpha$ [resp. $\alpha$], il saute ensuite
 dans $\C$ uniformément sur ses voisins et il allume [resp. éteint] à nouveau 
 la 
 lampe au point o\`u il se trouve avec probabilité 
 $1-\alpha$ [resp. $\alpha$].
 Ces  trois étapes sont indépendantes. Remarquons que si on "oublie"
 les configurations et que l'on regarde uniquement le premier argument de la
 marche $Z$, on retrouve la marche $X$, ainsi on peut écrire $Z_n=(X_n,f_n)$.
Par ailleurs, si nous fixons une trajectoire dans $\C$, les
états des différentes lampes
sont indépendants.
Enfin, la marche $Z$ admet des mesures réversibles, elles sont proportionnelles 
à:
$$m(a,f)= \nu (a) (\frac{1-\alpha}{\alpha})^{ \# \{  i ; f(i)=\bar{1} \}  }.$$
Nous verrons plus tard qu'il suffit d'étudier le cas $\alpha =
1/2$ pour obtenir la borne supérieure dans le théorème \ref{theoinitial}(la mesure $m$ se réduit alors à 
$m(a,f)=\nu(a)  $ ) mais pour la borne inférieure du théorème et pour d'autres
valeurs de $\alpha$, la mesure $m$ nous sera utile. 

Le lien, entre la marche $Z$ sur ce produit en couronne et notre problème
initial réside dans la propriété suivante. On note 
 $\tilde{\P}^{\omega}_o$ la loi de la marche $Z$  issue de $o$.
\begin{prop} 
\label{proba-Nn} On a, 
$\ \tilde{\P}^{\omega}_o (Z_{2n} =o) = \mathbb{E}_0^{\omega}
(\alpha^{N_{2n}} 1_{ \{ X_{2n}=0\} } )   .$
\end{prop}

\begin{proof}
\begin{eqnarray*}
\tilde{\P}^{\omega}_o (Z_{2n} =o) &=& \tilde{\P}^{\omega}_o \Bigl(  (X_{2n},
f_{2n} )=(0,f_0)   \Bigr) \\
                           &=&\underset{ \underset{ k_0=k_{2n}=0  }
{(k_0,k_1,...,k_{2n} )\in \Z^d} } {\sum}
                                  \tilde{\P}^{\omega}_o (X_0=k_0,
X_1=k_1,...,X_{2n}=k_{2n} \   \\
&\ & \hspace*{6.5cm}\text{et } f_{2n}=f_0)\\
                          &=& \underset{ \underset{ k_0=k_{2n}=0  }
{(k_0,k_1,...,k_{2n} )\in \Z^d} } {\sum}
                                 \tilde{\P}^{\omega}_o (X_0=k_0,
X_1=k_1,...,X_{2n}=k_{2n} ) \\
&\ & \hspace{2cm} \times \ \ \tilde{\P}^{\omega}_o( f_{2n}=f_0 |X_0=k_0,...,X_{2n}=k_{2n}) \\
                          &=&\underset{ \underset{ k_0=k_{2n}=0  }
{(k_0,k_1,...,k_{2n} )\in \Z^d} } {\sum}
                                                   \tilde{\P}^{\omega}_o (X_0=k_0,
X_1=k_1,...,X_{2n}=k_{2n} )\\
&\ & \hspace{6.5cm}\times \alpha^{ \# \{ k_0,...,k_{2n} \} }   \\
                         &=& \mathbb{E}_0^{\omega} (         \alpha^{N_{2n}} \
1_{       \{  X_{2n}=0 \}        }         ).
\end{eqnarray*}
\end{proof}

  \subsection{Remontée de l'isopérimétrie sur le produit en couronne.}
  \label{anna}
  On explique dans cette sous section, comment une inégalité
  isopérimétr-ique 
  sur le graphe $\C_n^g$ se transmet au produit en couronne $\C^g_n \wr \frac{\Z}{2\Z}$.
 Ce type de résultat est d\^u à A. Erschler. Nous présentons ici, une preuve
 légérement plus simple et plus détaillée. 
  Le résultat principal est le suivant:
\begin{prop}
\label{trucanna} Il existe des constantes universelles $C_1, C_2>0$ telles que,  
$$Fol_{\C_n^g \wr \frac{\Z}{2\Z} }^{\C^g \wr \frac{\Z}{2\Z}}
 (k)  \geq
e^{C_1 Fol_{\C_n^g}^{\C^g} (C_2k)}.$$
On peut prendre (cf preuve) $C_1= \log(2)/9$ et $C_2=1/1000$.
\end{prop}
  La preuve de cette propriété est assez  technique et découle de plusieurs
  lemmes. Pour alléger les notations, on notera, dans cette section, 
  $ \partial_{\wr} =\partial_{\C^g \wr \frac{\Z}{2\Z}} $. 
  Pour comprendre la preuve,  examinons d'o\`u proviennent les points du bord
  d'un ensemble  $U\subset V(\C^g_n\wr \frac{\Z}{2\Z})$.
  Il y en  a  de deux types: \\
- les points  du bord provenant du bord de la base $\C_n^g $, ie : les points de la 
forme  $\Bigl( (a,f);(b,f)\Bigr)$  avec  $(a,b)\in E(\C^g) $, \\
- les points du bord provenant du bord en 'configuration', ie : les 
points de la forme $\Bigl( (a,f);(a,g)\Bigr)$ avec $f=g$ sauf en  $a$.\\
Avant d'énoncer les lemmes préliminaires à la preuve de la propriété
\ref{trucanna},  introduisons quelques notions. Soit 
   $U \subset V(\C_n \wr \frac{\Z}{2\Z}) .$
   \begin{enumerate} [$\bullet$]
  \item A chaque $U$, on associe un graphe $K_U $ de la manière suivante :\\
- les points sont les configurations $f$ de l'ensemble 
$$\{f ;\  \exists a \in \C_n \ \ (a,f)\in U \}, $$ 
- deux points distincts $f$ et  $g$ sont reliés par une ar\^etes si
et seulement si :
$$ \exists a \in \C_n  \ \text{tel que }
  \left\lbrace
\begin{array}{l}
(a,g) \ \text{et} \ (a,f)  \in U,  \\
\text{et} \\
\forall x\neq a  \ f(x)=g(x).
\end{array}
\right.
$$
 ie: $f=g$ sauf en un point de $p(V)$, o\`u $p$ est la projection 
   $V(\C_n^g \wr \frac{\Z}{2\Z}) \rightarrow V(\C_n^g)=\C_n.$
\item Soit  $f\in K_U$, on dit que  $f$ est 
$b-satisfaisable$,  si il y a au moins $b$ ar\^etes attachées à  $f$  dans $K_U$,
$$ie: \ \ \ \#\{x\in p(V) ; \ \underset{x}{\dim} f =1      \} \geq b,$$  
o\`u
$$ \underset{x}{\dim} f = \#\{g ; \ (x,g)\in V \ \text{et}  \ \forall y \neq x \  f(y)=g(y) \} \in \{0,1\}.$$
Nous noterons  $S_V(b)$ l'ensemble des configurations $b-satisfaisables$, 
pour alléger les notations, la plupart du temps nous oublierons la dépendance 
en  $U$, et on notera donc  $S(b)$ cet ensemble.
\begin{rem}  \label{cardino}
Pour tout $U$, on a $$ |U| \geq 2 | E(K_U)|, $$
et si pour tout $(x,f) \in U$, on a $\underset{x}{\dim} f=1$ (
si $U$ n'a pas de points "isolés") alors,
$$ |U|=2   \# E(K_U ).  $$  
\end{rem}
\item Si $f$ n'est pas b-satisfaisable, on dira que $f$  est  $b-non\  satisfaisable$ et
on notera $NS(b) $ l' ensemble des configurations b-non satisfaisables.
\item Soit $K $  un graphe, on étend la notion de satisfaisabilité à $K$. Un
 point $x\in V(K)$ sera dit $b-satisfaisable$ [resp.
$b-non \  satisfaisable$] si il y a au moins [resp. strictement moins de] $b$ ar\^etes attachées à $x$.
\item une ar\^ete sera dite  $b-satisfaisable $  si elle relie deux
configurations  b-satisfaisables sinon on dira qu'elle est b-non satisfaisable. 
On notera  $S^e(b)$ ou $S^e_U(b)$ [resp. $NS^e(b)$] l'ensemble des 
ar\^etes  b-satisfaisables [resp. b-non satisfaisables]. 
\item Un point  $u =(x,f)\in U$  sera dit  $b-satisfaisable $
[resp. $b-non$ $satisfaisable $] si $f \in S(b)$ 
[resp. $ NS(b) $]. On utilisera la  notation $S^p(b)$ et  $NS^p(b)$  pour 
l'ensemble des points  de $U$  qui sont (ou ne sont pas) b-satisfaisables. \\ 
\item Un point $u=(x,f)\in U$ sera dit  $bon$ si 
 $ \underset{x}{\dim} f = 1$ sinon on dira qu'il est $mauvais $.
\end{enumerate}
Maintenant la preuve de  la propostion \ref{trucanna}  se décompose en trois étapes,  
on suppose  que  $\frac{|\partial_{\wr} U| }{| U| }
\leq \frac{1}{k}$, on prouve d'abord qu'il y alors peu de points
  $b-non\ satisfaisables $  (pour une  certaine valeur de  $b$), puis on extrait 
  un sous graphe tel que tous les points soient    $\frac{b}{3}-satisfaisables $
  , et enfin on en déduit une minoration de $|U|.$ \\
\\
Dans la proposition suivante, on note  $\phi$ la fonction $Fol_{\C_n^g}^{\C^g}$.

\begin{lem}
\label{neud} Soit  $U \subset V(\C_n^g \wr \frac{\Z}{2\Z}) $, tel que 
$\frac{|\partial_{\wr} U| } {| U| } \leq \frac{1}{1000k}$ alors 
  \begin{enumerate}[$(i)$]
\item  $  \frac{ |\{  u \in U ;\  u\ \text{mauvais}   \} | }{|U|}
\leq \frac{1}{1000k}, $
\item  $  \frac{|\{  u \in U ;\  u \in
NS^p(\phi(k)/3 ) \  \}| }{|V|}    \leq \frac{1}{500}. $
\end{enumerate}

\end{lem}
\begin{proof} $\ $ \\
Pour (i) on remarque qu'un point mauvais de $U$ nous donne un 
point du bord en 'configuration'  de $U$, ainsi:
\begin{eqnarray*}
|\partial_{\wr} U| \geq \#\{ u \in U;\ u \ \text{mauvais}  \}.
\end{eqnarray*}
Pour  (ii) posons:
\begin{eqnarray*} Neud &=&\{  u \in U ;\  u \in
NS^p(\phi(k)/3 ) \  \} \\
&=&\{ u =(x,f) \in U  ; f \in NS(\phi(k)/3 )\},
\end{eqnarray*}
et
 $$Neud(f) = \{ (x,f) ; (x,f)\in U  \}.$$  
Notons que l'on a   $p(Neud(f))=\{x; (x,f)\in U\}.$ Pour un ensemble  $F$  de 
 configurations, on pose  
   $Neud(F)= \underset{f\in F}{\cup} Neud(f) $. Remarquons que c'est une union
   disjointe.
   
   Soit maintenant  $f \in NS(\phi(k)/3)$,
   int\'eressons-nous à l'ensemble
    $p(Neud(f))$. Chaque point du bord de  $p(Neud(f))$ fournit un point du 
    bord  de $U$ en 'base'. Soit l'ensemble $p(Neud(f))$ poss\`ede  
     une part importante   de bord relativement à son volume, 
      et il fournit  alors une part de bord en 'base' dans $U$ de l'ordre de son
      volume. Ou bien, cet ensemble poss\`ede peu de bord relativement à son 
      volume, et dans ce cas, du fait que $f$ soit non satisfaisable, 
  cet ensemble fournit à $U$ du bord en 'configuration'. 
  
  Dans tous les cas, on obtient des points du  bord de $U$ mais les hypothèses sur $U$ limitent
  cet   apport. 
 On distingue donc deux cas. \\
 \\
\underline{Premier cas} : $f\in F_1 := \{  f\in NS(\phi(k)/3) ;   
\frac{| \partial_{\C^g} p(Neud(f)) |}{ | p(Neud(f)) |} > \frac{1}{k} \}. $\\

On fait correspondre à chaque point du bord de $p(Neud(f)) $ un point du bord de $U$
par l'application injective suivante:  \begin{eqnarray*}
 &\underset{f\in F_1}{\dot{\bigcup}}& \partial_{\C^g} p(Neud(f)) \longrightarrow 
  \partial_{\wr} U  \\
       &(x,y)& \longmapsto \Bigl(\left(x,f\right);\left(y,f\right)  \Bigr).
       \end{eqnarray*}
Donc, on a:
    \begin{equation}
 |\partial_{\wr} U  | \geq  \underset{f\in F_1}{\sum} | \partial_{\C^g}
p(Neud( f) )| \geq  \frac{1}{k} \underset{f\in F_1}{\sum} |
p(Neud(f))|
\geq  \frac{1}{k}|Neud (F_1)| . \label{cas1}
\end{equation}
\\
\underline{Second cas} : $f\in F_2 := \{  f\in NS(\phi(k)/3) ;  
 \frac{ |\partial_{\C^g} p(Neud(f)) |}{| p(Neud(f))| } \leq \frac{1}{k} \} .$\\
\\
Comme $f\in NS(\phi(k)/3 ) $, on a:
$$ | \{ x\in p(Neud(f)) ; \ \underset{x}{\dim} f =1 \} |< \frac{1}{3}\phi(k), $$
d'o\`u 
$$ |\{ x\in p(Neud(f)) ; \ \underset{x}{\dim} f =0 \} |\geq |Neud(f)| - \frac{1}{3}\phi(k). $$ \\
(on a utilisé le fait que  $|p(Neud(f)|=|Neud(f)|.$)\\
Or  $ f \in F_2$, donc par la définition de la fonction de F\o lner on a: 
 $$|Neud(f)| \geq \phi(k).$$
Ainsi, $$ | \{ x\in p(Neud(f) ; \ \underset{x}{\dim} f =0 \}| \geq \frac{2}{3} |Neud(f)| .$$

Mais pour tout  $f\in F_2$, si $y\in \{ x\in p(Neud(f) ; \ \underset{x}{\dim} f =0 \}$
 le point $(y,f)$  nous donne un point du  bord de $U$ en configuration 
 de manière injective, par l'application
 \begin{eqnarray*}
 &\underset{f\in F_2} {\dot{\cup} }& \{ (x,f)\in Neud(f) ; \ \underset{x}{\dim} f =0 \} \rightarrow \partial_{\wr } U \\
 &(x,f)& \longmapsto \Bigl( (x,f );(y,\bar{f}_x)  \Bigr),
 \end{eqnarray*}
 o\`u $\bar{f}_x$  est la configuration qui vaut $f$ en tout point différent de $x$
 et dont la valeur en $x$ est $\bar{0}$ [resp. $\bar{1}$] si $f(x)=\bar{1}$ 
 [resp. $f(x)=\bar{0}$]. 
 Donc, on obtient,
$$ |\partial_{\wr} U | \geq \frac{2}{3} \underset{f\in F_2}{\sum} |Neud(f)| \geq 
\frac{2}{3} |Neud(F_2)|.$$
D'o\`u l'on déduit que pour  $k\geq 2 $  \begin{equation}
  \label{cas2}
|\partial_{\wr} U | \geq \frac{1}{k} |Neud(F_2)|. \ \ \ 
\end{equation}
En additionnant  (\ref{cas1}) et  (\ref{cas2}) et en utilisant que 
 $\frac{|\partial_{\wr} U| } {|U| } < \frac{1}{1000k}$ , on a   : \\
$$ \frac{|Neud|}{|U|} < \frac{1}{500}.$$
\end{proof}

\begin{lem}
\label{soushypergraphe} Soit  $b >0$ et soit  $ K$ un   graphe. Supposons que 
 $E(K) \neq \emptyset$  et :\\
$$ \frac{  |NS_K^e(b) |}
{|E(K)| }<1/2.
$$\\
Alors il existe un sous graphe non vide  $ K'$  de $K$ tel que toutes les 
ar\^etes soient  $S_{K'}^p(b/3).$
\end{lem}
\begin{proof}
 On efface tous les sommets  $NS^p(b/3)$ et ainsi que toutes les ar\^etes 
 adjacentes. Il peut alors appara\^itre  de nouveaux sommets $NS^p(b/3)$ dans ce
 nouveau graphe. On les efface à nouveau ainsi que les ar\^etes adjacentes, et on
 itère ce processus  "d'effacement" tant qu'il y a des sommets $NS^p(b/3)$.
Prouvons que ce processus  se termine avant que le graphe soit vide. \\
Posons,
\begin{eqnarray*} C_1&=& |NS_K^e(b) | ,\\
C_2&= &| \{e \in S_K^e(b) ; e \text{ effacée à la fin du processus }  \} | \\
   C_0&=& | \{ e\in E(K); e \text{ effacée à la fin du processus } \}|.  
 \end{eqnarray*}
Si nous montrons que  $ C_2 \leq C_1$, le lemme sera démontré car
 $$C_0\leq C_1 +C_2 \leq 2 C_1<
 |E(K)|.$$ 
Ce qui signifiera qu'il reste  des points  non effacés.

Pour prouver que $C_2\leq C_1$, introduisons une orientation des ar\^etes. 
Si  $A$ et  $B$ sont deux points du graphe $K$ reliés par une ar\^ete, on oriente 
l'ar\^ete de $A$ vers $B$, si $A$ a été effacé avant $B$. Nous choisissons une
orientation arbitraire si $A$ et $B$ sont effacés simultanément ou si $A$ et $B$
sont non effacés.
Nous noterons  $\underset{\downarrow}{A}\ $ [resp. $\underset{\uparrow}{A}$]  l'ensemble
des ar\^etes orientées quittant $A$ [resp. arrivant en $A$]. 

\begin{souslem}
\label{soustrucchiant}
Soit  $b>0$ et   $A $ un point de  $K$, effacé à la fin du processus. 
Supposons que  $A$ soit initialement  $ S^p(b)$,  alors 
$$ |\underset{\downarrow}{A}|   \leq \frac{1}{2}| \underset{\uparrow}{A}|.   $$
\end{souslem}
%
%
%
%
\includegraphics*[width=11cm]{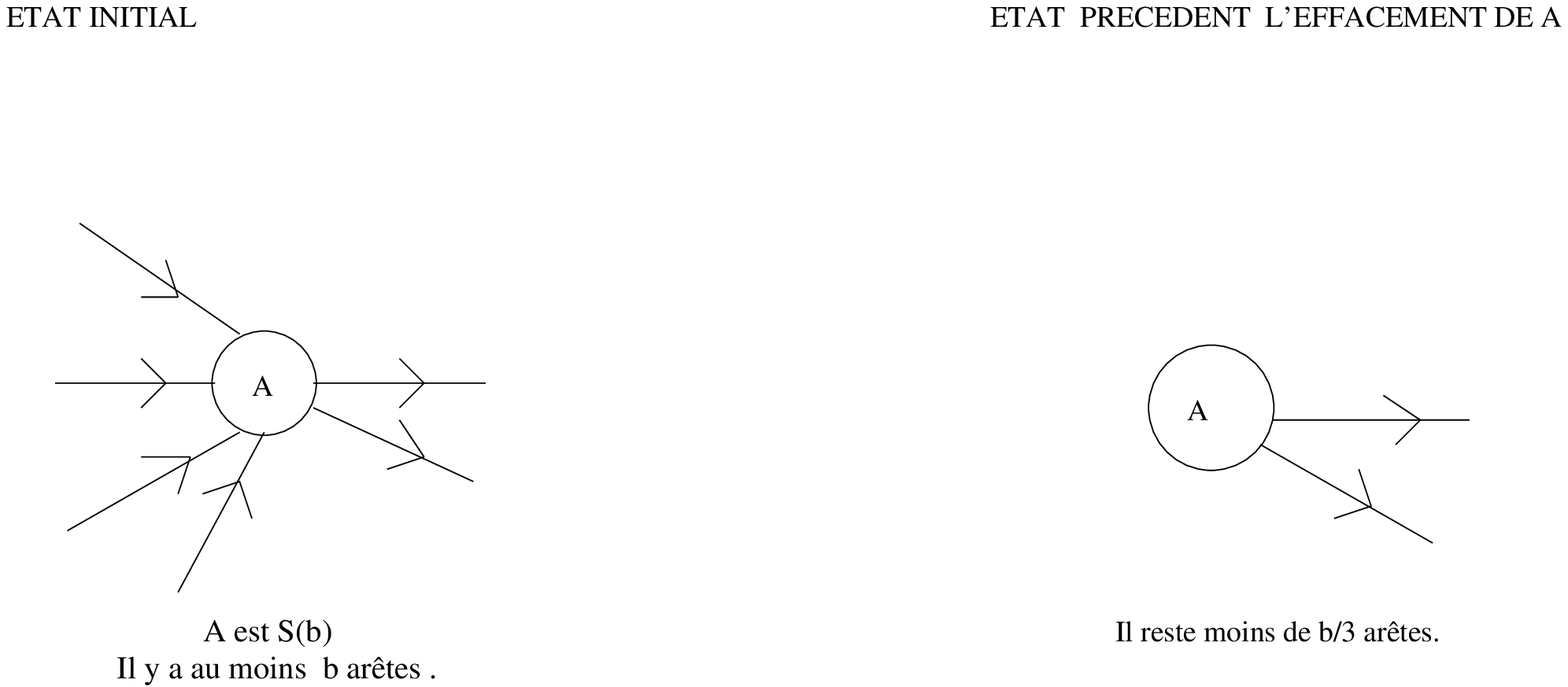} \\
$$\text{Figure a}$$
\begin{proof} La figure a illustre le début et la fin du processus au point $A$. 
On a,
 
$$ |  \underset{\downarrow}{A}|  \leq b/3  \ \ \text{et}\ \
 |\underset{\uparrow}{A} | \geq b-\frac{1}{3}b \geq \frac{2}{3}b.    $$
D'o\`u,  $$ |  \underset{\downarrow}{A} | \leq \frac{1}{2} |  \underset{\uparrow}{A} | . $$
\end{proof}
Pour finir la preuve, on pose:
\begin{eqnarray*}
D_1&=&\{ \text{points effacés à l'étape 1}  \}, \\
D_i&=&\{ \text{points  initialement $ S^p(b)$, effacés à l'étape i} \} \  \text{pour } \ i\geq 2 \\
F_i&=& \{ \text{ar\^etes entre  \ $D_i$ et $ D_{i-1}$ }  \}, \\
 F_i'&=& \{ \text{ar\^etes quittant  \ } D_{i-1}  \}.     
 \end{eqnarray*}
Ces ensembles sont schématisés sur la figure b.\\


\includegraphics*[width=11cm]{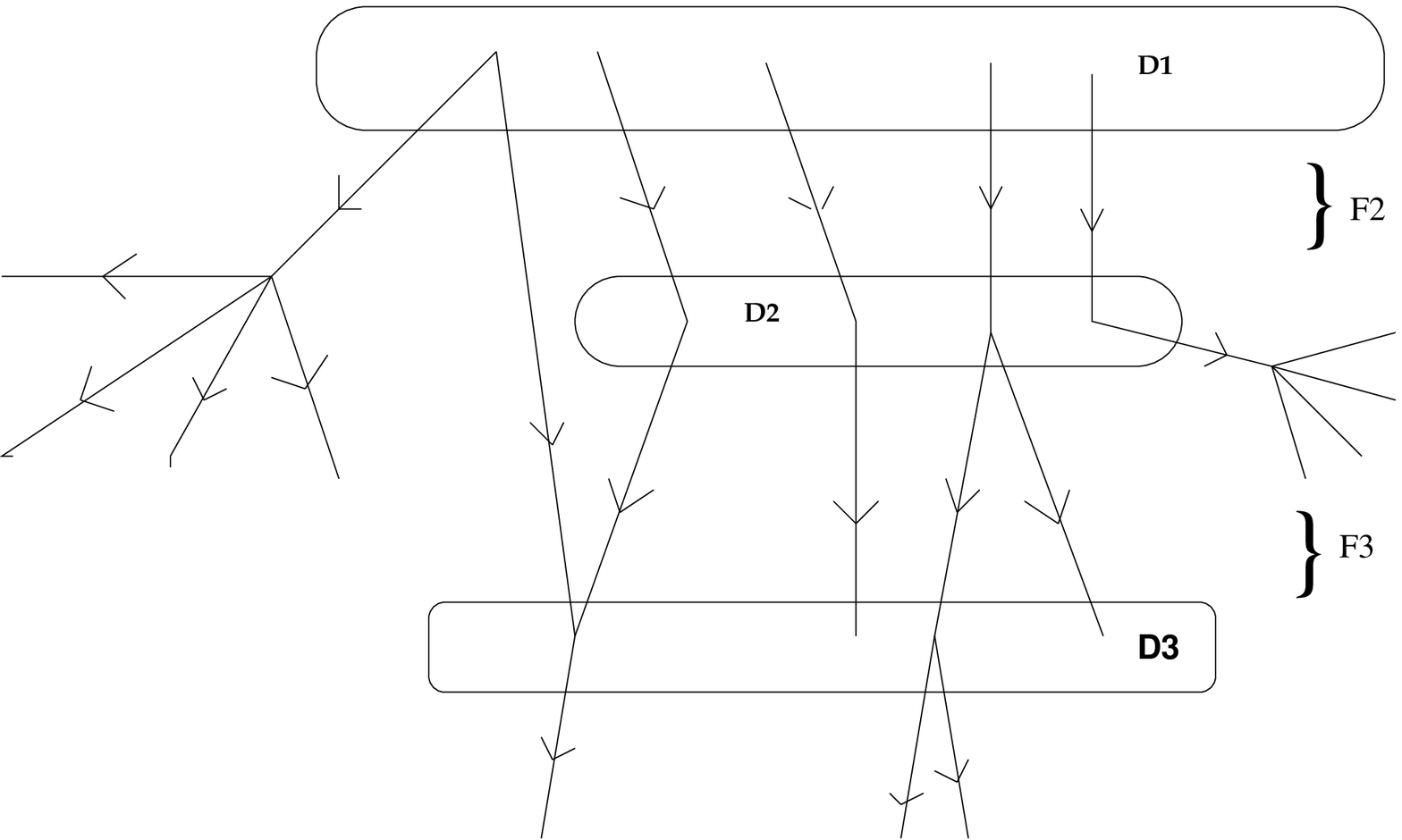}
$$\text{Figure b}.$$
Notons que  $F_i \subset F_i'$  et que par ailleurs les ar\^etes de   $F_i'$ 
sont effacées à l'étape $i-1$.
En appliquant  le sous lemme \ref{soustrucchiant},  en chaque point de 
 $D_i$, on obtient:
$$\forall i \geq 2 \ \  |F_{i+1}'|  \leq \frac{1}{2} | F_i| . $$
D'o\`u, $$ \ \  |F_{i+1}'  | \leq \frac{1}{2^{i-1}}| F_2|.  $$
(On a utilisé le fait que  $|F_i|\leq|F_i'|.$)
Ainsi,  \begin{eqnarray} \label{ca}  | \underset{i\geq 3}{\cup} F_i' |   \leq
(\underset{i\geq 1} {\sum} \frac{1}{2^i}) \ | F_2|.   \end{eqnarray}

Maintenant, les ar\^etes de $F_2$ sont  $NS_K^e(b)$  car si elles \'etaient 
 $S_K^e(b)$, elles relieraient deux points  $S_K^p(b)$ et en particulier,
 les  points de $D_1$  seraient $S^p_K(b)$,  donc $S^p_K(b/3)$  donc non 
 effacés à l'étape 1. Il s'en suit que,
 \begin{eqnarray} \label{F2}| F_2| \leq | NS^e_K(b)|=C_1 .
 \end{eqnarray}
Par ailleurs, toute ar\^ete  $S^e(b)$ effacée, est dans un certain $F_i'$ 
avec  $i\geq 3$, donc
\begin{eqnarray} \label{Fi'} \ \ \  \ \ C_2=\# \{e \in S^e(b) ; e \text{ effacée à la fin du processus }  \}      
\leq  |\underset{i\geq 3}{\cup} F_i'| .
\end{eqnarray} 
Des inégalités (\ref{Fi'}), (\ref{F2}) et (\ref{ca}), on déduit que 
$C_2\leq C_1$, ce qui termine la preuve.
\end{proof}

\begin{lem}
\label{choixsuccessif} Soit $Y >0$ et soit  $\A$ un ensemble non vide de 
configurations, tel que pour toute configuration de  $\A$, il existe au moins 
$Y$ points o\`u l'on peut changer la valeur de $f$ tout en restant dans $\A$,
alors  : $|\A| \geq 2^Y.\ $  i.e :\\
 $ (\forall f \in \A \ \exists a_1,a_2,...,a_Y \in \C_n \text{ tels que }  \ \ 
 \bar{f}_{a_i} \in \A) 
 \Longrightarrow  |\A|\geq 2^Y,$
\\
o\`u $\bar{f}_{a_i}$ est définie à partir de  $f$ par: 
\begin{eqnarray*} 
\bar{f}_{a_i} (x)=
\begin{cases}
 \;  f(x) & \text{si $x \neq a_i $} \\
  \; \bar{1}-f(x) & \text{sinon . }
\end{cases}
\end{eqnarray*}
\end{lem}
\begin{proof}
On raisonne par récurrence sur  $Y$.\\
-Si  $Y=1$  la propriété est satisfaite.\\
-Supposons $Y\geq 2$ et considérons un point  $x$  (dans la
base ) tel qu'il existe  $f  , g\in A $ satisfaisant   $f(x)=0 $ et  $g(x)=1 $, et posons 
$\A_0 = \{h\in \A ; h(x)=0   \} $ et  $\A_1 = \{h\in \A ; h(x)=1   \}. $ 
 $\A_0$ et $\A_1$  sont non vides. On a 
$\A=A_1 \cup A_0$ et cette union est disjointe. Par ailleurs,  $A_0 , \ \A_1$  v\'erifient 
l'hypothèse de 
récurrence avec la constante  $Y-1$. 
D'o\`u, $$|\A| =| \A_0| + |\A_1| \geq 2.2^{Y-1} =2^{Y+1}.$$
\end{proof}
A l'aide des trois lemmes précédents, on peut maintenant
démontrer la propriété
\ref{trucanna}.\\
\\
\textbf{Preuve de la propriété  \ref{trucanna}:} \\
 Soit $U\subset  \C_n^g \wr \frac{\Z}{2\Z}$ vérifiant 
 $\frac{|\partial_{\wr} U|}{|U|}\leq \frac{1}{1000k}$.
 Procédons en 5 \'etapes.
\begin{enumerate}[1.] 
\item Soit $U'=U-\{u\in U;\ u \  \text{mauvais}\}$.
Posons   $\tilde{K}=K_{U'}  $, on a $ E(K_U) =E(K_{U'})$ ($\tilde{K}$ est le 
sous graphe de $K_U$ qui ne contient que des points attachés à une ar\^ete 
de $K_U$).
   \item $V(\tilde{K})$ (et donc aussi  $E(\tilde{K} )$  par construction) est
   non vide, puisque par le  (i) du lemme \ref{neud} :
 $$|V(\tilde{K}) |  \geq
(1-\frac{1}{1000k})|U|.$$
\item
 On a successivement:
\begin{eqnarray*}
 \# \{e\in E(\tilde{K})\cap NS^e(\phi(k)/3)\} )
  &=&
\#\{e\in E(K_U ) \cap NS^e(\phi(k)/3)\} \\
  &\leq& \frac{1}{2} \#\{ u\in U ; NS^p\Bigl(\phi(k) /3\Bigr)  \}
\\
&\ &   \hspace{0cm}\mathrm{ (remarque \ \ref{cardino})}\\
&\leq &  \frac{1}{1000}|U  |  \\
&\ & \hspace{0cm} \text{ (par le lemme  \ref{neud}} \ (ii) )\\
  &\leq&  \frac{1}{ 1000-\frac{1}{k} } \ \#\{u\in U;\ u \ \text{bon} \} \\
&\ &  \hspace{0cm} \text{ (par \ le  \ lemme \  \ref{neud} }\   (i) )\\
 &= &\  \frac{2}{ 1000-\frac{1}{k} }   | E(\tilde{K})|  \\
 &\leq &\  \theta   \  \ |E( \tilde{K}) | ,   \\
\end{eqnarray*}
avec  $\theta= \frac{2}{ 999 }<  \frac{1}{2}$.\\
\item Gr\^ace aux points 2 et 3, on peut appliquer le lemme \ref{soushypergraphe} à $\tilde{K}$ 
et déduire qu'il existe un sous graphe $K'$ 
de $\tilde{K}$  o\`u tous les points sont  $S^p_{K'}(\phi(k)/9)$.\\
 \item On déduit du  lemme \ref{choixsuccessif} que,  
 $$|U| \geq 2^{ \phi(k)/9}.$$
Ainsi on a prouvé que pour tout ensemble $U\subset V(\C_n^g\wr\frac{\Z}{2\Z})$,
$$ \frac{|\partial_{\wr} U|}{|U|}\leq\frac{1}{1000k}  \Rightarrow 
|U|\geq 2^{ Fol_{\C_n^g}^{\C^g}(k)/9},$$
ce qui implique,
$$  Fol_{\C_n^g\wr\frac{\Z}{2\Z}}^{\C^g\wr \frac{\Z}{2\Z}} (k)
 \geq 2^{ Fol_{\C_n^g}^{\C^g}(k/1000)/9},$$
 et fournit des valeurs numériques des constantes $C_1$ et $C_2$.
\end{enumerate}

\section{Isopérimétrie sur  un amas de  percolation }
\label{section is}
Après avoir bri\`evement discuté d'une première inégalité isopérimétrique, 
on explique quel type d'inégalité isopérimétrique est nécessaire pour prouver 
le théorème \ref{theoinitial}. Enfin le reste de cette section  est consacrée
à la preuve de  la propriété \ref{is de base}.
\subsection{Point de départ.}
Dans  \cite{pierre}, l'inégalité suivante est prouvée.
\begin{prop} \label{is1}
Il existe une constante $\beta=\beta(p,d) >0 $ telle que  $Q$ p.s sur l'ensemble
$ |\C|=+\infty$ on ait,
$$ \exists n_{\omega} \ \ \forall n\geq n_{\omega } \ \ 
\underset{A\subset \C_n , |A|\leq |\C_n|/2}{inf} \
\frac{|\partial_{\C_n^g}A|} {|A|^{1-\frac{1}{\epsilon}}} \geq
\frac{\beta}{
n^{1-\frac{d}{\epsilon} }  },$$ 
o\`u  $\epsilon=\epsilon(n) = d+2d\frac{\log\log (n)}{\log(n)}$.
\end{prop}

Dans cette inégalité,  on peut remplacer le bord dans $\C_n^g$ par le 
bord dans $\C^g$, et enlever la condition sur le volume de
$A$, en remarquant les deux points
suivants.
\begin{enumerate} [1.]
\item  Pour tout $A\subset C_n$,  $   \partial_{\C^g} A = \partial_{\C_{n+k}^g}A$ pour tout  $ k\geq 1$.
 \item Puis pour enlever la condition sur le volume,  rappelons la propriété
 relative à la
 croissance du volume de $\C_n$ ( voir  Appendix B dans \cite{pierre})
 \begin{prop} \label{C^n} Il existe $\rho > 0$ tel que $ Q$  p.s  sur l'ensemble  
 $|\C|=+\infty$  pour $n$ assez grand, $$|\C_n| \geq \rho n^d.$$
   \end{prop}
  Ainsi, par exemple il existe  $ c'>0$  tel que  $Q$ p.s pour $n$ assez grand,
  pour tout $A\subset C_n$ on ait, $$|A| \leq |\C_n| \leq |\C_{n+c'n}| /2 .$$ 
En effet   $|\C_n|\leq n^d\leq  \frac{\rho}{2} (n+k)^d \leq \frac{1}{2} |\C_{n+k}|$
est réalisée dès que l'on prend  $k\approx n$. \end{enumerate}
On déduit donc de la propriété \ref{is1}, (en modifiant la constante $\beta$ par
une constante multiplicative) que:
  $$\exists \beta>0 \  Q \ \text{ p.s pour $n$ assez grand,} \   
   \underset{A\subset \C^n }{inf} \
\frac{|\partial_{\C^g} A|} {|A|^{1-\frac{1}{\epsilon}}} \geq
\frac{\beta}{ n^{1-\frac{d}{\epsilon} }  }.$$ 
 Cette manière de compter le bord nous donne une majoration de 
 la probabilité de retour de la marche $Z$ tuée quand elle
 sort de la bo\^ite
 $\B_n=[-n,n]^d$.  Si on applique alors la m\^eme  démarche 
  que l'on va utiliser à la
 section 4, à partir de cette inégalité isopérimétrique, nous trouvons  une
  majoration
 de la probabilité de retour de  $Z$ (et donc de la  transformée de Laplace 
 du nombre de points visités)  en  $$e^{-c  \frac
{
  t^ {\frac{d}{d+2} }
 }
 {     {\log(t)}^{c'}    }
  }.$$

Ceci  n'est pas la borne supérieure attendue. La principale partie du travail
consiste donc à supprimer le terme logarithmique.
 Pour cela, l'idée est de prouver une nouvelle inégalité isopérimétrique, 
 qui est similaire à celle que l'on a dans le graphe  $\L^d=(\Z^d,E^d)$  
    pour les ensembles de 'gros' volume.    Plus précisement, la propriété
    \ref{is de base} peut s'interpréter ainsi,  pour  $A\subset \C_n $, 
   si $|A|$ est grand,  on a
    $\frac{ |\partial_{\C^g} A |}{ |A|^{1-\frac{1}{d} }} \geq
C$, et si  $|A|$ est petit, on peut dire que $|\partial_{\C^g} A| \geq 1 $.

Les deux sections suivantes sont consacrées à prouver la propriété \ref{is de base} pour tout $p>p_c$, par
des techniques de renormalisation.  Dans la section
\ref{other}, nous prouvons pour des valeurs 
de $p$ proches de 1, une   inégalité isopérimétrique 
modifiée, correspondant 
à un événement croissant. Puis dans la section \ref{renorm},
à l'aide du  théorème 2.1 de   \cite{antal}, on déduit la
propriété \ref{is de base}.

\subsection{Une autre inégalité isopérimétrique}\label{other}

Soit   $\G_n $ l'ensemble des points de $\B_n$ attachés à une ar\^ete
ouverte. Soit alors  ${\L}_n$ la plus grosse composante connexe 
de  $\omega$ dans  $\B_n$. 
On note $\overline{\L_n}$ la composante connexe de $\L_n$
dans le graphe $\omega $. 
  $\overline{\L_n^g}$  est le graphe ayant
  $\overline{\L_n}$
  comme ensemble de points et  l'ensemble des ar\^etes est défini par: 
  $E(\overline{\L_n^g} )=\{ (x,y); \ \omega(x,y)=1, \
  x,y\in\overline{\L_n} \}. $
   
  Posons pour  $A \subset \B_n$,   $n(A)$ le nombre de composantes connexes
  de  $\B_n-\L_n$  qui contiennent au moins une composante connexe de  $A$. 
  Notons que si  $A$ est  connexe $n(A) =0 \ \text{ou} \ 1$
  suivant que $A\subset \L_n$ ou non.\\
Remarquons  également que si $A\subset\B_n$,  
 $|\partial_{\bar{\L_n^g} } A| =
 |\{ (x,y)\in E^d ; x \in A\cap \overline{\L_n}  ,\ 
 y\in B_{n+1}-A, \text{ et } \omega(x,y)=1
       \}|. $

\begin{prop} \label{propISnew}
Il existe  $p_0 < 1$ tel que pour tout  $p \in [p_0;1]$:
\begin{eqnarray} \nonumber
\quad \quad  \exists \beta>0 \ \ \forall c>0 \ \ 
Q \text{ p.s }  \exists n_{\omega,c} \ \ \ 
\forall n\geq n_{\omega,c} \ \ \ \\
\label{ISnew}
 \forall A \subset \B_n \  connexe  \ \ \ \frac{ n(A) + |\partial_{\bar{\L_n^g}} A |}{ f_c(|A|)} \geq \beta.
 \end{eqnarray}
\end{prop}
Avant de prouver cette propriété, remarquons les deux  faits suivants:\\
-d'abord l'événement  défini par l'équation (\ref{ISnew}) est croissant. En effet, l'ajout d'une ar\^ete à
 $\omega$, ne fait pas diminuer  $|\partial_{\bar{\L_n^g}} A |. $
Et si  $n(A)$ diminue de 1, cela signifie qu'il y avait une composante connexe
de  $A$ qui n'intersectait pas  $\L_n$ et qui maintenant intersecte $\L_n$, 
donc le bord  $|\partial_{\bar{\L_n^g}} A|$a  augmenté de 1. Finalement, 
la somme  $n(A) + |\partial_{\bar{\L_n^g}}
A |$ ne décro\^it pas.\\
C'est cette raison qui motive l'introduction d'une telle
inégalité, afin de permettre l'utilisation des techniques de
renormalisation. \\
-par ailleurs, pour  $n$ assez grand, 
$\C_n=\L_n$ (voir la première partie de la preuve de
\ref{propISnew}),  
   si l'on prend alors $A\subset \C_n$ connexe,  $ n(A)=0$ et on retrouve 
l'inégalité  \ref{is de base}  pour $p$ proche de 1, comme
conséquence de la proposition \ref{propISnew}.

\begin{proof}  $ \ $\\
Commen\c{c}ons par prouver  que pour $n$ assez grand
$\L_n=\C_n$, ce qui nous sera utile dans la suite de la
preuve.\\
Supposons $\L_n\neq \C_n$. On a  immédiatement par la
propriété \ref{C^n} que:
\begin{eqnarray} \label{L_n}
|\L_n|\geq |\C_n|\geq \rho n^d.
\end{eqnarray}
Soient $C_i$ les composantes connexes de  $\B_n -\L_n$ dans le graphe $\B_n^g$. 
Pour tout  $i$ on a:
 \begin{eqnarray*}
\begin{cases} \; \partial_{\B_n^g} C_i  \text{ est 
(*)-connexe\footnotemark[2] } \\
\; (\text{car }C_i  \text{ et } \B_n-C_i \text{ sont connexes } )\\
\;\text{et}\\ 
\; \frac{ \underset{e\in  \partial_{\B_n^g} C_i } {    \sum}   1_{  \{ w(e)=1 \} }  }
{ |\partial_{\B_n^g} C_i|}=0.
\end{cases}
\end{eqnarray*}
\footnotetext[2]{un ensemble $M$ d'ar\^etes est dit
(*)-connexe si pour tout $e,e'$ de $M$, il existe une suite
d'ar\^etes $e_1,...,e_m $ de $M$ telle que $e_1=e$ et
$e_m=e'$ et $e_i \overset { (*)  }{\sim} e_{i+1}$ où      $(x,y)\overset { (*)  }{\sim}(x',y')$ si 
$ \max |x_i-x_i'| \leq 1 \  \text{et} \   \max
|y_i-y_i'|\leq 1.   $}
Pour $\beta >0$, posons,
\begin{eqnarray} 
\A_n'= \{       
\exists F \subset E(\B_n^g) ; \    
 \frac{ \underset{e\in  F } {    \sum}   1_{  \{ w(e)=1 \} }  }
{ |F|} \leq \beta \ ,
\end{eqnarray} 
$$ \hspace*{6.cm} |F|\geq (\log n)^2 , F \  
\mathrm{(*)-connexe}
\}.$$
On a:
\begin{eqnarray*}   Q(\A_n' )\leq 
\underset{  \underset{ |F| \geq(\log n)^2 \text{ et } F \
(*)-\text{connexe}} {F}           } {\sum}
  Q(\frac{ \underset{e\in  F } {    \sum}   1_{  \{ w(e)=1 \} }   }
{ |F|} \leq \beta) 
\end{eqnarray*}
Du fait que les variables aléatoires $(\omega(e) )_{e}$
suivent des  des lois de Bernouilli(p) et sont
indépendantes,  par l'inégalité de
Bienaymée Tchebitchef, on  a:\\ 
\begin{eqnarray}\label{bienayme}\forall \lambda\geq 0 \ \ \  Q(\frac{ \underset{e\in  F } {    \sum}   1_{  \{ w(e)=1 \} }   }
{ |F|} \leq \beta) \leq  e^{\lambda \beta |F|}(pe^{-\lambda} + 1-p)^{|F|} .
\end{eqnarray}
D'où, pour tout $\lambda \geq 0$ on déduit:
\begin{eqnarray} \label{q}  Q(\A_n' )\leq 
\underset{  k \geq(\log n)^2            } {\sum}
 e^{\lambda \beta k }(pe^{-\lambda} + 1-p)^{k}  
   \times
\end{eqnarray}
$$ 
 \hspace{4cm}|\{F \subset E(\B_n^g );\ 
F\ \mathrm{(*)-connexe},\ |F|= k \}|.      
$$
Mais  $ |\{F \subset E(\B_n^g );\ 
F\ \mathrm{(*)-connexe},\ |F|= k \}|
\leq (2n+1)^d e^{ak},$  o\`u $a$ est
une constante qui dépend seulement de $ d$ (voir
\cite{sinai}).
Ainsi, le terme général de la série du membre de droite de l'équation (\ref{q})
est majoré par:
$$ (2n+1)^d e^{k[a+\lambda\beta+\log(pe^{-\lambda} +1-p)]  }.$$
 Si $p$ est proche de $1$, et si on choisit $\beta$ assez petit (avec 
 $\beta< 1$), il existe  $\lambda >0$  tel que 
  $ a + \lambda \beta + \log(pe^{-\lambda} +1-p) = - \xi < 0. \ $
Donc pour $n$ assez grand, 
$$ Q(\A_n' )\leq C\ (2n+1)^d e^{- \xi  (\log n)^ 2   }\leq  e^{- \frac{\xi}{2} 
 (\log n )^ 2   }  
$$
Cette dernière expression est sommable en $n$, donc par le lemme de Borel
Cantelli, on déduit: 
\begin{eqnarray*}
\ \forall \omega \ \exists n_{\omega} \ , \ \forall n\geq n_{\omega} \ 
\ \omega \in \A_n'^c .
\end{eqnarray*}
 Ainsi, $Q$ p.s pour $n$ assez grand 
$\partial_{\B_n^g}C_i \in A_n'^c$.   
\begin{eqnarray} \label{bordC_i} ie: \ \ |\partial_{\B_n^g} C_i| \leq   ( \log n)^2 .
\end{eqnarray}
Or $$|C_i|\leq |B_n| - |\L_n|.$$
Puis de (\ref{L_n}), on déduit:
$$|C_i|\leq c' n^d \text{ avec } c'<1.$$
Puisque $C_i$ a un volume  plus petit
qu'une fraction du volume de $\B_n$, on
peut donc lui appliquer l'inégalité
 isopérimétrique restreinte à une boite  $\Z^d$, pour en déduire
\begin{eqnarray} \label{volumeC_i}
|C_i| \leq  C_d \log(n)^{\frac{2d}{d-1}}.
\end{eqnarray}
Or, si $\L_n\neq \C_n,$ il existe $i$ tel que $\C_n\subset
C_i$. Ce qui est impossible compte tenu 
de  (\ref{volumeC_i}) et de la
propriété \ref{C^n}.
Donc pour $n$ assez grand, $\L_n=\C_n$.\\

Fixons-nous à présent un $c>0$.\\
Si $|A| < c n^{\gamma}$ , (\ref{ISnew}) est satisfaite puisque 
$n(A) + |\partial_{\bar{\L_n^g} } A | \geq 1.$ \\
Supposons donc maintenant 
$|A| \geq c n^{\gamma}$. Par la majoration (\ref{volumeC_i} ) vérifiée par les
volume des $C_i$, nécessairement $A$ intersecte $\L_n$.
De ce fait, on a $n(A)=0$.
Nous 
voulons donc  prouver que
\begin{eqnarray} \label{onveut}
\frac{|\partial_{\bar{\L_n^g}} A|}{|A|^{1-\frac{1}{d}}}  \geq \beta.
\end{eqnarray} 
Procédons en 5 étapes.
\begin{enumerate} [1.]
\item L'isopérimétrie classique dans  $\Z^d$ fournit, 
\begin{eqnarray} \label{isA}
\frac{|\partial_{\L^d} A|}{|A|^{1-\frac{1}{d}}}  \geq C_d.
\end{eqnarray}
\item On voudrait remplacer  $\partial_{\L^d} A$  par 
 $\partial_{\bar{\L_n^g}} A$. Mais si on applique un argument de contour à
 $A$, il peut exister des ar\^etes de  $\partial_{\L^d} A$  qui appartiennent 
 à $\omega-\overline{\L_n^g}$. Nous allons donc appliquer un argument de contour à l'ensemble 
 $A$ auquel nous avons rajouté les
 "morceaux" de $\omega-\overline{\L_n^g}$ connexes à $A$. \\
   Du fait que p.s, il n'existe qu'une seule composante 
 connexe infini dans $\omega$  ( qui est $\C_n=\L_n$ pour  $n$ assez grand), 
 on peut trouver un  $k=k_{n,\omega}$ tel que tous les chemins 
de  $\omega -  \C_n^g$ commen\c{c}ant dans $\B_n$, se terminent 
avant d'atteindre une ar\^ete  de  $\partial_{\L^d} \B_{n+k}$.\\

$$\includegraphics[width=7cm]{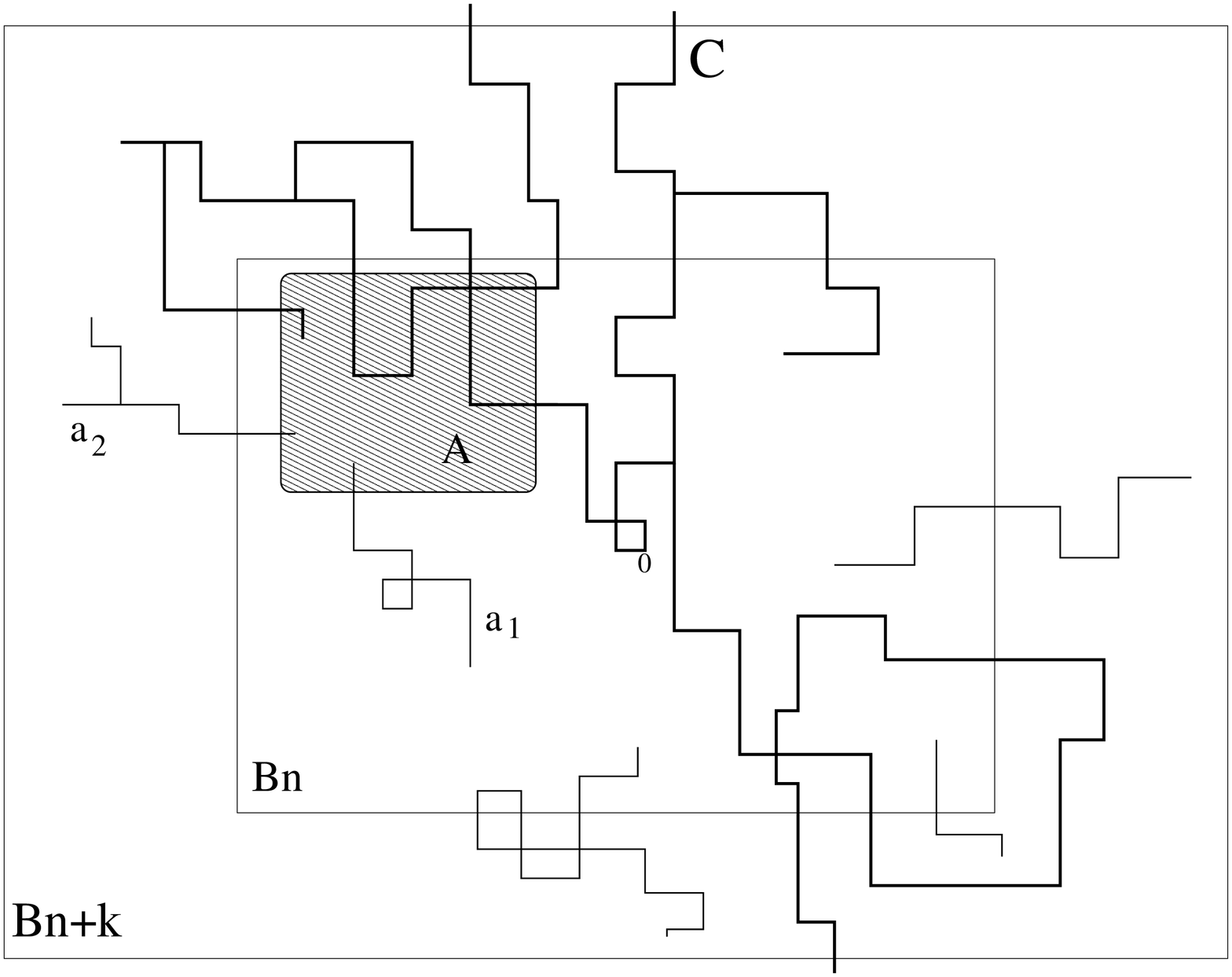} $$
$$\text{figure d.}$$

Soient $(a_i)_{i=1...m}$ les chemins de  $\omega -
\overline{\L_n^g}$  qui débutent
dans  $ A$  (voir figure d). Tous ces chemins sont contenus 
dans $\B_{n+k}$. Ajoutons ces chemins à $A$, et posons 
 $$A''=A \cup a_1...\cup a_m .$$
 \item $A''$ vérifie les propriétes suivantes. $A''$ est connexe et on peut 
 supposer que $\Z^d-A''$ est connexe. En effet, si ce n'est pas 
 le cas, on remplit les trous, en prenant $B$  le
 complémentaire dans $\Z^d$ de la composante
 connexe de $\Z^d-A''$ (dans le graphe $\L^d$) qui contient les points à l'infini, et
 une inégalité isopérimétrique pour l'ensemble $B$ se
 transmet à $A''$, puisque $|B|\geq |A''|$ et 
 $|\partial_{\L^d}  B| \leq | \partial_{\L^d}  A''|$ .\\
 Ainsi $\partial_{\L^d} A''$
 est (*)-connexe .
 
 
 Par ailleurs, on a   \begin{eqnarray}
 \label{bordchiant}|\partial_{\omega} A''| =|
 \partial_{\bar{\L}_n^g} A|,
 \end{eqnarray}
($\omega$  désigne dans l'égalité précédente le graphe issu du processus de percolation) 
  et   l'inégalité  isopérimétrique classique dans $\Z^d$ donne :  
 $$|\partial_{\L^d} A''| \geq  C_d |A''|^{1-\frac{1}{d}} 
 \geq C_d (cn^{\gamma})^{1-\frac{1}{d } }. $$
 \item  Appliquons un argument de contour à $A''$. Pour $\beta'>0$ on pose:
  \begin{eqnarray*}
  \hspace*{0cm} \A''_n= \{ \omega; \        
\exists F \subset E^d , \ |F|\geq C_d
(cn^{\gamma})^{1-\frac{1}{d } }   , \ 
 \frac{ \underset{e\in  F } {    \sum}   1_{  \{ w(e)=1 \} }  }
{ |F|} \leq \beta', \\
 \hspace*{2.5cm}  F \  \mathrm{(*)-connexe}
,\  F \cap E(\B_n^g)\neq \emptyset \}. 
\end{eqnarray*}
Comme précédemment, à l'aide de (\ref{bienayme}) on a:
$$ Q(\A''_n )\leq 
\underset{j \geq C_d (cn^{\gamma})^{1-\frac{1}{d } }   } {\sum}
e^{\lambda \beta' j}(pe^{-\lambda} + 1-p)^{j} \times h(n,j),$$
o\`u $ h(n,j)=  \underset{k\geq 0}{\sup} |\{ H \subset
E^d \ ,
H\ \mathrm{(*)-connexe} \ 
 |H|= j \ , \\
\hspace*{7cm} \text{et} \ H \cap E(\B_n^g)\neq \emptyset  \}|.      $\\
\\
Mais $ |\{ H \subset E(\B_{n+k_{n,\omega}}^g) \ ,
H\ \mathrm{(*)-connexe} \ , \  |H|= j \ , \ H \cap E(\B_n)\neq \emptyset  \}|
     $ est inférieur à  $(2n+1)^d e^{aj}$ pour tout $k\geq 0$ (avec $a$
     indépendant de $k$). Donc,  on obtient finalement que si  $p$ 
     est proche de $1$, il existe une valeur  $\beta' >0$, et 
    un $\delta>0$ tels que pour $n$ assez grand:
$$ Q(\A''_n ) \leq \  e^{- \delta n^{\gamma (1-\frac{1}{d })}  } . $$ 
Puis par le lemme de Borel-Cantelli,  on d\'eduit que  $Q$ p.s pour $n$ assez grand,
$\partial_{\L^d} A'' \not\in \A_n''$.
 $$ie: \ | \partial_{\omega} A'' | \geq  \beta' |\partial_{\L^d} A'' | . $$
 \item On peut alors terminer la preuve:  par (\ref{bordchiant}) on a,
 $$| \partial_{\bar{\L_n^g}} A | \geq  \beta' |\partial_{\L^d} A'' | . $$
En appliquant alors, l'inégalité isopérimétrique dans $\Z^d$ à $A''$
et par le fait que $A\subset A''$, on déduit qu'il existe  $\beta>0$ tel que
$Q $ p.s, pour $n$ assez grand, 
$$ | \partial_{\bar{\L_n^g}} A | \geq  \beta |A|^{1-\frac{1}{d} } .$$
  \end{enumerate}
 \end{proof}

\subsection{Renormalisation} \label{renorm}

Nous sommes maintenant en mesure de montrer  la propriété \ref{is de base} pour tout 
 $p>p_c$. Reprenons les notations et résultats de
 \cite{antal}.
Soit  $p>p_c$ et   $N \in \N$. On recouvre $\Z^d$  par une union disjointe 
de boites de taille  $N$ telles que 
$\Z^d=\underset{i \in \Z^d}{\cup} B_i$ o\`u $B_i$  est la boite centrée
en $(2N+1)i$.
Soit $B_i'$ la boite de m\^eme centre que $B_i$ mais de taille $\frac{5}{4}N$. 

Introduisons (voir \cite{antal}), 
 $\mathcal{E}(N)=\{  (ke_1,(k+1)e_1 ) ; k=0,...,[\sqrt{N}] \} $ et  
$\mathcal{E}_i=\tau_{i(2N+1)} \mathcal{E}(N), $  o\`u $\tau_b$  représente 
le  décalage  dans  $\Z^d$  avec $b \in \Z^d$.
On dit qu'un amas $K$ contenu dans $B'$
est traversant pour $B\subset B'$, si  dans les $d$ directions il existe un
chemin ouvert contenu dans $K\cap B$ reliant la face droite de $B$ à la face
gauche de $B$.

 Considérons les deux événements
suivants:
\begin{eqnarray*}
R_i^N :&=&\{  \text{il existe  un amas traversant } K \text{ dans } B_i' 
  \text{ pour } B_i', \text{ tout che-} \\  
  &\ & \text{ min ouvert contenu  dans} 
\  B_i' 
 \text{ de longueur plus grande que}   \\
 &\ & \ N/10  \text{ est connecté à }  
  K \text{ dans }   B_i', 
       \text{ et }   K \text{ est
       traversant pour} \\ 
&\ & \text{ toute sous-boite }      B \subset B_i' \text{ de taille supérieure à } 
        N/10 \},
	 \end{eqnarray*}
\\
$S_i^N:\ \  = \ \ \{ \text{ Il y a au moins une ar\^ete ouverte dans } \mathcal{E}_i 
\}.$ \\

Nous écrirons  $K(B_i')$  quand cela est nécessaire pour indiquer 
que nous considérons l'amas traversant $K$ de la boite $B_i'$.

On définit alors  $\phi_N : \omega \longrightarrow \{ 0,1\}^{\Z^d} $ 
telle que  
$(\phi_N \ \omega)(x)= 1_{R_x^N \cap S_x^N }(\omega) ,$ et on dit  
qu'une boite $B_x$ est  $bonne$ si $(\phi_N \ \omega)(x)=1.$ (ie: 
si les deux événements sont réalisés) sinon on dit que la boite est  $mauvaise$.
Cela définit un processus (dépendant) de percolation sur les  sites
du réseau renormalisé, voir  \cite{antal} pour plus de détails.
\\
\\
{\bf Preuve de la proposition \ref{is de base}.}\\
On a prouvé la propriété \ref{propISnew}
pour un paramètre de percolation sur 
les ar\^etes assez grand. Puisque l'ensemble
des amas qui vérifie l' inégalité (\ref{ISnew}) est un événement
croissant,  cette inégalité  est donc aussi vérifiée pour un paramétre de 
percolation sur les sites  assez grand (voir \cite{liggett}). 

Soit  $A \subset \C_n$ connexe,  notons 
 $\tilde{A}=\{ i \in \Z^d ; A \cap B_i \neq \emptyset \}$.  $\tilde{A}$ est 
 connexe et remarquons que l'on  a:
 $$ |\tilde{A}| \geq \frac{ |A|}{ (2N+1)^d}.$$
 On utilise la notation $\ \tilde{} \ $  pour les quantités 
 définies dans le processus renormalisé.
 
Pour un bon choix de $N$,  par la proposition 2.1 de  \cite{antal}, 
on peut appliquer la  proposition \ref{propISnew} pour  $\tilde{A}$ dans 
le processus renormalisé,  donc  (\ref{ISnew} )  est vérifiée pour  $\tilde{A}$:
il existe $\beta>0$ tel que pour tout $c>0$, $Q$ p.s pour $n$ assez grand on
ait, 

\begin{eqnarray} \label{isatilde} \frac{ n(\tilde{A}) +
|\partial_{\bar{\tilde{\L_n^g }}} \tilde{A} |}
{ f_c(|\tilde{A} |)} \geq \beta. 
  \end{eqnarray}

On distingue alors deux cas.
\begin{enumerate} [1.]
\item
Si $\tilde{n} (\tilde{A})= 1$, alors $\tilde{A}$ est contenu dans une
des composantes connexes de  $\tilde{\B_n}-\tilde{\L_n}$. Ainsi,
 $\partial_{   \bar{  \tilde{\L}   }_n^g     } \tilde{A}=\emptyset, $  et donc 
  (\ref{isatilde}) donne $1 \geq \beta f_c( |\tilde{A}|)$. Du fait que 
   $\partial_{\C^g} A \neq \emptyset $, on déduit : 
$$ |\partial_{\C^g} A| \geq \beta f_c( |\tilde{A}|) \geq \beta 
f_c(  \frac{|A|}{ (2N+1)^d}) 
\geq  \frac{\beta}{(2N+1)^{d-1 } }  f_{ c(2N+1)^d} (|A|). $$
 qui est l'inégalité de la proposition \ref{is de base} 
 (pour tout $c>0$ et avec une constante  $\beta$ plus petite).  \
\item Si  $\tilde{n}(\tilde{A})=0$,  alors  par (\ref{isatilde}) on a, 
 $$\frac{  |\partial_{\bar{\tilde{\L} }_n^g} \tilde{A} |}{ f(|\tilde{A} |)}
 \geq \beta. $$   
Soit  $(i,i') \in  \partial_{\bar{\tilde{\L} }_n^g} \tilde{A} $ alors
 $B_i $ et  $B_{i'}$ sont des bonnes boites et  $A$
intersecte $B_i$ mais pas  $B_{i'}$. Nous allons montrer que chaque couple
ainsi choisi donne au moins une ar\^ete de  $\partial_{\C^g } A$.
En effet, du fait que les deux boites $B_i$ et $B_{i'}$ sont bonnes, on peut trouver
   $x \in B_i \cap A $ et  $ y \in B_{i'} \cap (\C-A)$
reliés par un chemin ouvert dans  $B_i \cup B_{i'}$, et sur ce chemin se trouve
une ar\^ete  de $\partial_{\C^g } A$.

[ Expliquons rapidement comment trouver  $x$ et  $y$. On prend  $x\in A$ 
donc  $x\in B_i \cap \C $.  $\C$ est infini donc il existe un chemin dans 
 $\C$ reliant  $x$ à l'extérieur de  $B'_i$. Comme ce chemin a une longueur
 supérieure à  $N/4$,  il est connecté à l'amas  $K(B_i')$  dans $B'_i$ 
 (donc par connexité $K(B_i') \subset \C $).
Prenons maintenant un chemin ouvert contenu dans  $K(B_i') \cap B'_i$ joignant
la face gauche à la face droite  de la boite  $B'_i$ dans la direction
de  $B_{i'}$  et soit  $y $ sur ce chemin et dans  $ B_{i'}$.
Toujours par connexité, on en déduit finalement  l'existence d'un chemin dans  $\C$ 
 reliant  $x$ à $y$. ] 
 
Il y a au plus  $2^d$ ar\^etes  $ (i,i')$   
qui peuvent donner la m\^eme ar\^etes dans  $\partial_{\C^g } A$ (toutes les
ar\^etes $(i,i')$ avec $i'$ voisin de $i$ dans $\Z^d$ ). 
Donc finalement, on a, pour tout $c>0$:
 $$    |\partial_{\C^g} A | \geq 
 \frac{1}{2^d} |\partial_{\bar{\tilde{\L} }_n^g}  \tilde{A}|
  \geq \frac{ \beta}{2^d } f_c(|\tilde{A }| ) \geq  
 \frac{\beta}{2^d (2N+1)^{d-1} }  f_{ c(2N+1)^d} (|A|)  . $$
\end{enumerate}
Ce qui termine la preuve de la propriété  \ref{is de base}, pour tout 
 $p \geq p_c.$\\

\section {Preuve de la borne supérieure}
Par les résultats de la section \ref{section is} et \ref{anna}, on prouve une
minoration de la fonction de F\o lner sur $\C_n^g\wr \frac{\Z}{2\Z}$, puis par des
outils analytiques (type inégalité de Nash) on obtient une borne supérieure de la probabilité
de retour de $Z$ puis de la transformée de Laplace du nombre
de points visités. 
\subsection{Isopérimétrie sur $\C_n^g \wr \frac{\Z}{2\Z}$. }
La propriété \ref{is de base}, donne la minoration suivante pour la fonction de 
F\o lner de $\C_n^g$ relative à $\C^g$. 
\begin{prop}
\label{ISbase}
Soit $\gamma>0$, il existe $\beta >0$ tel que pour tout
$c>0,$  $\ Q \ p.s $  sur 
l'ensemble  $|\C|=+\infty$, pour $n$ assez grand, on ait:
\begin{eqnarray}  \label{folnebaseprop }
 Fol_{\C_n^g}^{\C^g}(k) \geq
\begin{cases}
 \; k & \text{si $k < c n^{\gamma}$} \\
  \; (\beta k)^d & \text{si  $k\geq   c n^{\gamma}$.}
\end{cases}
\end{eqnarray}\\
\end{prop}

\begin{proof}
On proc\`ede  en deux étapes.
\begin{enumerate}[1.]
\item D'abord il est facile de voir que l'on peut restreindre 
le minimum définissant la fonction de
F\o lner à des ensembles connexes,
\begin{eqnarray*}
Fol_{\C_n^g}^{\C^g}(k) &=& \min \{  |A| ; A\subset \C_n \ \text{et } \ 
\frac{|\partial_{\C^g} A| }{| A| } \leq
 \frac{1}{k}\}\\ &=&
\min \{  |A| ; A\subset \C_n \ A \ \text{connexe et } \ 
\frac{|\partial_{\C^g} A| }{| A| } \leq
 \frac{1}{k}    \}. 
    \end{eqnarray*}
En effet, si $A$ n'est pas connexe, notons  $(A_i)_{i=1..t}$ ses 
composantes connexes. Au moins une satisfait la condition 
$\frac{|\partial_{\C^g} A_i| }{| A_i| } <
 \frac{1}{k},  $   et on a  $ |A_i| \leq |A| .$ 
\item Maintenant, par la proposition \ref{is de base}, il existe 
$\beta >0$ tel que pour tout $c>0, \ Q \ p.s$ sur  
$|\C|=+\infty $ pour $n$ assez grand, on a 
 \begin{eqnarray*} 
  \frac{|\partial_{\C^g} A | }{f_c(|A|) } \geq
\beta \ \ \text{ pour tout ensemble } A\subset\C_n \text{ connexe},
\end{eqnarray*}
Donc, pour tout $A$ connexe de $\C_n$, on a,
\begin{eqnarray} \label{impl}\frac{|\partial_{\C^g} A | }{|A| } \leq \frac{1}{k}   
\Longrightarrow  \frac{ |A|}{f_c(|A|)}  \geq \beta k .
\end{eqnarray}
Posons alors  $\G (x)=\frac{x}{f_c(x)}$  et définissons $\G^{-1}$ par: 
$$ \G^{-1} (y) = \inf \{ x \  ; \G(x)=y \} =
\begin{cases}
\;y  &\text{si } x<cn^{\gamma},\\
\; y^d &\text{si } x\geq cn^{\gamma}. 
\end{cases}$$
On déduit de (\ref{impl})  $$|A| \geq \G^{-1} (\beta k)  .$$
\begin{eqnarray} \label{folenfin} ie:\ Fol_{\C_n^g}^{\C^g} (k) \geq \G^{-1} (\beta k ) . \end{eqnarray}
\end{enumerate}
 \end{proof}
Gr\^ace à la proposition  (\ref{trucanna}) et à l'inégalité (\ref{folenfin}), 
on déduit imm\'ediate-\\ment qu'il existe $C >0$, tels que pour tout $c>0, \ Q$ p.s  pour $n$ assez grand,    
    \begin{eqnarray}
    \label{FolW}    \ \
     Fol_{\C_n^g \wr \frac{\Z}{2\Z}}^{\C^g\wr\frac{\Z}{2\Z}} (k)  \geq
  \begin{cases}
  \;e^ {Ck} & \text{si }   k < cn^{\gamma}, \\
 \;e^{Ck^d}   & \text{si }   k \geq  cn^{\gamma} .
\end{cases}
\end{eqnarray}
On ne peut déduire directement une majoration des noyaux de la marche $Z$, à
partir de cette minoration car les ar\^etes du graphe
 $\C^g\wr \frac{\Z}{2\Z}$ ne sont pas adaptées aux sauts de $Z$. On construit
 donc les graphes suivants.

Soit $G$ un graphe, on note $G\wr\wr \frac{\Z}{2\Z}$ le graphe tel que,
$$\begin{cases}
\;V(G\wr\wr \frac{\Z}{2\Z} )=V(G\wr \frac{\Z}{2\Z}),\\ 
\;E(G\wr\wr \frac{\Z}{2\Z})=\{\Bigl( (x,f)(y,   f_{x,i \atop y,j}    )\Bigr) ; 
\  (x,y)\in E(G) \ \text{et }  \ (i,j)\in\{ 0,1\} \ \\
\hspace{8.7cm} \text{ et supp} (f) \subset G \ \}.
\end{cases} $$
Soit $\delta $ un point imaginaire,  on pose
$\overline{\C_n}=\C_n\cup\{\delta\}$ et
$\overline{\C_n^g}$ le graphe ayant $\overline{\C_n}$
comme ensemble de points et l'ensemble des ar\^etes est défini par
$E(\overline{\C_n^g})=\{ (x,y);\ \omega(x,y)=1,\
x,y\in\C_n \}\cup\{ (x,\delta); \ x\in\C_n,\ \exists
z\in V(\omega) \ (x,z)\in\partial_{\omega} \C_n\}$.
 Notons alors
$$ W_n= \overline{\C^g_n} \wr \frac{\Z}{2\Z}$$
et 
$$W_n'= \overline{\C^g_n} \wr\wr \frac{\Z}{2\Z}.$$
On prouve:
\begin{prop} Il  existe une  constante
$C>0$ telle que que pour tout $c>0, \ Q$
  p.s  sur $|\C|=+\infty$ et pour  $n$ assez grand,
\begin{eqnarray}
\label{fol-final}
Fol_{\C_n^g \wr \wr \frac{\Z}{2\Z}}^{ \C^g \wr \wr \frac{\Z}{2\Z}} (k) 
\geq F(k):=F_{C,c}(k)=
\begin{cases}
\;  e^ {Ck} & \text{si }  k <  c n^{\gamma}, \\
 \;e^{C k^d} & \text{si } k \geq c n^{\gamma} .\\
\end{cases}
\end{eqnarray}
\end{prop}
\begin{proof}
On proc\`ede en 4 étapes.
\begin{enumerate} [1.]
\item Soient $\epsilon$  et $\epsilon '$ les formes 
de  Dirichlet respectives  de $W_n$ et  $W_n'$ définies   
 pour tout 
$f,g : V(\overline{\C_n^g} \wr \frac{\Z}{2\Z}) \longrightarrow \R $
 par: $$\epsilon (f,g)= \underset{(x,y)\in E(W_n) }{\sum} (f(x)-f(y))(g(x)-g(y) ),$$
et
$$\epsilon' (f,g)= \underset{(x,y)\in E(W_n' )  }{\sum} 
(f(x)-f(y))(g(x)-g(y) ).$$

\item Les graphes  $W_n$ et $W_n'$ sont quasi isométriques
par l'application identité sur $V(W_n)=V(W_n')$.
En effet, soient $d$ et $d'$ les distances respectives sur ces
graphes, pour tout $u,v$ de $V(W_n) $ on a:
$$ \frac{1}{3} d'(u,v) \leq  d(u,v) \leq  3d'(u,v). $$
   Ainsi (voir I.3 dans
\cite{woess}), il existe 
$C_1, C_2 >0$ (indépendant de  $n$) tels que, pour tout 
$f : V(\overline{\C_n^g} \wr \frac{\Z}{2\Z}) \longrightarrow \R, $
\begin{eqnarray}\label{drich}  C_1 \epsilon (f,f) \leq \epsilon '(f,f)
 \leq  C_2 \epsilon (f,f).\end{eqnarray}

\item Soit $U\subset V(\overline{\C_n^g} \wr \frac{\Z}{2\Z})$, prenons $f=1_{\{    U \}}$
dans (\ref{drich}), on obtient
\begin{eqnarray}
\label{bordquasi-isom}
 C_1 |\partial_{ W_n  } U|  \leq   
 |\partial'_{W_n'      } U| \leq  
  C_2|\partial_{ W_n      } U|
\end{eqnarray} 

\item Or  on a $ |\partial_{ W_n  } U| = |\partial_{\C^g\wr \frac{\Z}{2\Z}} U |$
et $ |\partial_{ W_n ' } U| = |\partial_{\C^g\wr\wr \frac{\Z}{2\Z}} U |$.
Finalement, on déduit de  (\ref{FolW}) et (\ref{bordquasi-isom}),  que l'on
peut remplacer $\partial_{\C^g\wr \frac{\Z}{2\Z}    } $  par 
 $\partial_{\C^g\wr \wr\frac{\Z}{2\Z}}$ dans  (\ref{FolW}), 
 cela n'affecte que les constantes  $c $ et $C$. 
\end{enumerate}
\end{proof}
\subsection{Borne supérieure de $\tilde{\P}^{\omega}_o( Z_{2n}=o )$}
 On suppose   pour cette section $\alpha = 1/2$. Gr\^ace à (\ref{fol-final}),
 nous obtenons une majoration de la probabilité de retour de la marche 
 aléatoire $Z$. 
 \begin{prop} Supposons $\alpha=1/2$.
\label{resultproba}
Il existe une constante $ c>0 $ tel que $Q$  p.s   sur $|\C|=+\infty$, et 
 pour  $n$ assez grand, $$\tilde{\P}_o^{\omega}(  Z_n=o) 
  \leq  e^ {-c n^{\frac{d}{d+2} }} . $$
\end{prop}

\begin{proof}
Posons $\tau_n = \inf\{ s \geq 0 \ ; \ Z_s \not\in \C_n \wr \frac{\Z}{2\Z}  \}.$ 
Ainsi, on peut écrire,
\begin{eqnarray}
\label{casseen2}\tilde{\P}^{\omega}_o(Z_{2n}=o) =
\tilde{\P}^{\omega}_o(Z_{2n}=o  \ \mathrm{et} \ \tau_n \leq n)  
\end{eqnarray}
\begin{eqnarray*} 
\hspace{5.4cm} + \tilde{\P}^{\omega}_o(Z_{2n}=o  \ \mathrm{et} \ 
\tau_n > n   ).\end{eqnarray*}

-Le premier terme du membre de droite est nul, puisque la marche ne peut 
sortir de la "boite" $V(\C_n^g \wr\wr \frac{\Z}{2\Z})$ avant le temps $ n$.

-Le deuxième terme est égal à  $\P_o(\bar{Z}_{2n}^n = o),$ o\`u 
 $(\bar{Z}_i^n)_i$ désigne la marche aléatoire sur  $W_n'$, co\"incidant avec 
 $Z$ mais tuée quand elle sort de $W_n'$. $(\bar{Z}_i^n)_i$  est réversible pour la mesure  
 $m$ restreinte 
 à   $V(\C_n^g \wr\wr \frac{\Z}{2\Z})$ ($m$ est définie à la sous section  \ref{RW}).\\
   Par les méthodes 
 développées dans \cite{coulhon}), pour obtenir une borne supérieure 
  de $\P_o(\bar{Z}_{2n}^n = o)$, posons pour $ U \subset
  V(\C_n^g \wr\wr \frac{\Z}{2\Z})$ : 
   \begin{eqnarray*}
 \begin{cases}
 \;|\partial_{\C^g \wr\wr \frac{\Z}{2\Z} } U |_m =
  \underset{u_1,u_2}{\sum} m(u_1)\tilde{p}(u_1,u_2) \ 
  1_{\{ (u_1,u_2)\in \partial_{  \C^g \wr\wr \frac{\Z}{2\Z}  
   }U  \}} \\
  \;\mathrm{\ et \ }\\
 \;|U|_m = \underset{u\in U}{\sum}  m(u).
 \end{cases}
 \end{eqnarray*}
 Comme on suppose  $\alpha = 1/2$, la mesure  $m$ se réduit à 
   $m(a,f)= \nu(a)$ (qui est bornée car la valence du graphe  $\C^g$ est 
   comprise entre $1$ et $2d$)
  et on déduit donc  à l'aide de l'inégalité (\ref{fol-final}) que dans
  ce cas,   il existe une  constantes
  $C>0$ telle que pour tout $c>0,\ Q$ p.s sur 
 $|\C|=+\infty$, 
 pour $n$ assez grand, 

\begin{eqnarray}
\label{fol-end}
\ \ \ \ \ \min\{ |U|_m;\ U\subset V(\C_n^g \wr\wr \frac{\Z}{2\Z}) \ 
\text{et } 
 \frac{|\partial_{  \C^g \wr\wr \frac{\Z}{2\Z}   } U|_m}{|U|_m} 
 \leq \frac{1}{k}  \}
\geq F_{C,c}(k).
\end{eqnarray}
\\
  Prenons $c=1$. $F:=F_{C,1}$ est positive et croissante, et son inverse est:
  \begin{eqnarray} \nonumber 
 F^{-1} (y)&=& \inf\{ x \geq 0 \ ; \ F(x) \geq y \} \\
\label{F^{-1}} &=&\begin{cases}
 \;  \frac{1}{C} \log (y) & \text{si $y < e^{Cn^{\gamma}} $} \\
  \; c n^{\gamma}   & \text{si $ e^{Cn^{\gamma}} \leq y  \leq  e^{C n^{d \gamma}}    $}\\
  \;  (\frac{1}{C} \log(y) )^{1/d}   & \text{ si $ e^{C n^{d\gamma}}   <y,$}\\
\end{cases}
\end{eqnarray}

Comme  $F^{-1}$ est croissante, (\ref{fol-end})  implique  une
inégalité de  Nash  (voir par exemple 14.1 in \cite{woess}) pour la forme de 
Dirichlet associée à  $\bar{Z}^n$, puis une majoration du noyau: 
$$\text{pour tout } i\geq 0, \ \   \P_0(\bar{Z}_{2i}^n =0) \leq 2a(i),$$ 
o\`u $ a $ est la solution de l'\'equation différentielle suivante,
$$\left\lbrace
\begin{array}{l}
 a'=-\frac{a}{8(F^{-1} (4/a))^2},  \\
 a(0)=1. \\
\end{array}
\right.$$
A l'aide de l'expression de $F^{-1}$ par  (\ref{F^{-1}}) et  en résolvant l'équation différentielle
dans chaque cas, on obtient qu'il existe des constantes $c_i$ (avec 
$c_3,c_9 \ \text{et} \ c_{13}  >0$) telles que:
\begin{enumerate} [$\bullet$]
\item Pour $0 \leq t \leq c_1 n^{3 \gamma} + c_2 $,
$$  a(t) = e^{- (c_3  t + c_4)^{1/3} } .$$
\item Pour $ c_1 n^{3 \gamma} +c_2 \leq  t \leq   
 c_5 n^{(d+2) \gamma} + c_6 n^{3 \gamma}+c_7 n^{2 \gamma} +c_8  ,$
 $$     a(t) =    e^{ -c_9 t/n^{2 \gamma} +c_{10} n^{ \gamma} + c_{11}/n^{2
 \gamma} + c_{12}}. $$
 \item Pour  $  c_5 n^{(d+2) \gamma} + c_6 n^{3 \gamma}+c_7 n^{2 \gamma} +c_8    
  \leq t,$
$$ a(t) = e^{ -( c_{13 } t +c_{14}   n^{(d+2)\gamma} + c_{15} n^{3 \gamma} +c_{16}  n^{2 \gamma} +
   c_{17} )^{d/d+2}}  .   $$
\end{enumerate}
On choisit   $\gamma < \frac{1}{d+2} \ $, ainsi
 on déduit qu'il existe une constante $c>0$ telle que
    pour $t=2n$ assez grand,
$$ \tilde{\P}^{\omega}_o(\bar{Z}_{2n}^n = o)    \leq  e^ {-c n^{\frac{d}{d+2} }}.$$
Finalement par (\ref{casseen2}) et comme 
$n\mapsto \tilde{\P}^{\omega}_o(Z_{n} = o)$ est une fonction décroissante, on obtient  l'existence d'une
constante $c>0$ telle que $Q$ p.s sur $|\C|=+\infty$  pour $n$ assez grand,
$$ \tilde{\P}^{\omega}_o(Z_{n} = o)  \leq  e^ {-c n^{\frac{d}{d+2} }}  .$$
\end{proof}
\subsection{Conclusion: borne supérieure pour la transformée de Laplace.} 
Nous pouvons alors démontrer la borne supérieure de la propriété
\ref{theoinitial}.
\begin{prop}
\label{result} Pour tout $\alpha >0$, il existe une constante 
$ c(d,\alpha,p)>0 $ telle que $ Q \ p.s$ sur l'ensemble  $ |\C|=+ \infty$, 
$$\exists n_{\omega} \  \ 
\forall n \geq n_{\omega} \ \  \ 
\mathbb{E}_0^{\omega}(\alpha^{N_n})  \leq  e^ {-c(d,\alpha) n^{\frac{d}{d+2} }}.
$$
\end{prop} 
\begin{proof}
La résultat découle des 3 faits suivants.
\begin{enumerate} [1.]
\item D'abord, il suffit de prouver le résultat pour une valeur  $0<\alpha_0 <1$ 
de $\alpha$. 
En effet, supposons que 
$\mathbb{E}_0^{\omega}(\alpha_0^{N_n})  \leq  
e^ {-c(d,\alpha_0) n^{\frac{d}{d+2} }}$, alors

-si  $\alpha \leq \alpha_0 $,  il est clair que  
$\mathbb{E}_0^{\omega}(\alpha^{N_n})  \leq 
\mathbb{E}_0^{\omega}( \alpha_0^{N_n}) \leq
 e^ {-c(d,\alpha_0) n^{\frac{d}{d+2} }},$
 
-si $  \alpha_0 < \alpha < 1$, nous pouvons trouver  $\lambda >0$ tel que 
 $ \alpha={\alpha_0 }^{\lambda}$, avec $0< \lambda < 1$. Puis, on écrit,
\begin{eqnarray*}
\mathbb{E}_0^{\omega}(\alpha^{N_n}) &=& \mathbb{E}_0^{\omega}( \    [ \alpha_0^{  N_n} ]^{ \lambda}\ ) \\
&\leq& [\mathbb{E}_0^{\omega}( \     \alpha_0^{  N_n} \ )]^{ \lambda}  \ \
 \text{  \ (inégalité de Jensen)}
\\
&\leq& e^{-\lambda \ c(d,\alpha_0)  \ n^{\frac{d}{d+2} }}.  
\end{eqnarray*}
Ce qui donne le résultat pour $\alpha$ en prenant $c(d,\alpha)=\lambda
c(d,\alpha_0)= ( \frac{\log(\alpha)}{\log(\alpha_0)})c(d,\alpha_0)$.
\item La deuxième étape est résumée dans le lemme suivant (voir  \cite{SP}).
\begin{lem} \label{4.5}
Il existe une valeur $ \alpha_1  > 0$ et une constante $c_0'>0$  
telles  que,
 $$\mathbb{E}_0^{\omega}((1/2)^{N_{2n}} 1_{\{ X_{2n}=0   \}}  )  \geq 
c_0' \ \mathbb{E}_0^{\omega}(\alpha_1^{N_n}).$$
\end{lem}
\begin{proof} Pour "retirer" la condition  $\{ X_{2n} = 0\}$, l'idée est 
de couper les chemins à l'instant $n$ et d'utiliser  la réversibilité
entre l'instant  $n$ et $2n$. 
On utilise en particulier le fait suivant:  
\begin{eqnarray} \label{faitsalof}
 [\P_0^{\omega}( N_n=m)]^2  \leq   2d (2m+1)^d 
 \P_0^{\omega}( N_{2n} \leq 2m ;X_{2n}=0).
 \end{eqnarray}
Ecrivons en effet,
\begin{eqnarray*}
[\P_0^{\omega}( N_n=m)]^2 &=&   \Bigl(   \underset{ h\in B_m(\C)}{\sum}    
\P_0^{\omega} (N_n=m ; X_n=h)
\Bigr) ^2  \\
&=& 
   \Bigl(   \underset{ h\in B_m(\C)}{\sum}    
 \sqrt{\nu(h)}  \times 1/\sqrt{\nu(h)}  \times \\
  &\ & \hspace{4.4cm} \P_0^{\omega}(N_n=m ; X_n=h)
\Bigr) ^2    \\
&\leq&  \nu(B_m(\C)) \underset{ h\in B_m(\C) }{\sum} (1/\nu(h))   \P_0^{\omega}(N_n=m ; X_n=h)^2\\
& \ &\ \hspace{0cm} 
\text{ (par l'inégalité de  Cauchy-Schwarz) } \\
&\leq& 2d (2m+1)^d   
 \underset{ h\in B_m(\C) }{\sum}    \P_0^{\omega}(N_n=m ; X_n=h) \times \\
 &\ &  \hspace{3.5cm} \P_h^{\omega}(N_n=m ; X_n=0)(1/\nu(0)) \\
& \ &\ \hspace{0cm}  
\text{ (par   réversibilité) } \\
&\leq& 2d (2m+1)^d   
 \underset{ h\in B_m(\C) }{\sum}     \P_0^{\omega}(N_n=m ; X_n=h) \times \\
&\ & \hspace{3.3cm}  \P_0^{\omega}(N_n^{2n}=m ; X_n=h ; X_{2n}=0)\\
 & \ &\     \hspace{0cm}    ( \text{o\`u } N_n^{2n} = \#\{  X_n,...,X_{2n} \} )\\
  &\leq& 2d (2m+1)^d \ \P_0^{\omega}( N_{2n} \leq 2m;X_{2n}=0). 
\end{eqnarray*}
Calculons maintenant, $\E_0^{\omega}((1/2)^{N_{2n}} 1_{\{ X_{2n}=0   \}}  )$. On
a successivement:
\begin{eqnarray*}
\E_0^{\omega}((1/2)^{N_{2n}} 1_{\{ X_{2n}=0   \}}  ) 
&=& \underset{  m \geq 1}{\sum}   (1/2)^m \  \P_0^{\omega}( N_{2n}=m ; X_{2n}=0 )\\
&=& 1/2 \underset{  m \geq 1}{\sum}   (1/2)^m \ \P_0^{\omega}( N_{2n} \leq m ;
X_{2n}=0 ).\\
\end{eqnarray*}
Car $  \{ N_{2n}=m\}=\{N_{2n} \leq m\}-\{N_{2n}\leq m-1\}.$ Puis,

\begin{eqnarray*}
 \E_0^{\omega}((1/2)^{N_{2n}} 1_{\{ X_{2n}=0   \}}  ) 
 &\geq&  1/2 \underset{  m \geq 1}{\sum}   (1/4)^m \P_0^{\omega}( N_{2n} \leq 2m ; X_{2n}=0 )\\
&\ & \hspace*{0cm}(\text{par le fait (\ref{faitsalof})} \  )\\
&\geq&  \underset{  m \geq 1}{\sum} \frac{1}{4d (2m+1)^d}(1/4)^m [\P_0^{\omega}( N_n=m)]^2 \\
&\geq&   c_0 \underset{  m \geq 1}{\sum} (1/5)^m [\P_0^{\omega}( N_n=m)]^2 \\
&\geq&   c_0   \Bigl( \underset{  m \geq 1}{\sum} (1/4)^m
\Bigr)^{-1} \times \\
& \ & \hspace{2.8cm}\Bigl(
\underset{  m \geq 1}{\sum} (\frac{1}{2\sqrt{5}})^m\P_0^{\omega}( N_n=m) 
\Bigr)\\
&\ & \hspace*{0cm} (\text{par l'inégalité de Cauchy-Schwarz})\\
&\geq&   c_0' \ \E_0^{\omega} [ (\frac{1}{2\sqrt{5}})^{   N_n}] .\\
\end{eqnarray*}
ce qui prouve le lemme \ref{4.5}.
\end{proof}
\item On peut alors conclure. Pour $\alpha=1/2$, la propriété 
\ref{resultproba} donne
 $$ \P_o^{\omega}(  Z_n=o) 
  \leq  e^ {-c n^{\frac{d}{d+2} }}    .$$
  Par ailleurs, par la propriété \ref{proba-Nn}
  on a   
$$  \E_0^{\omega}((1/2)^{N_{2n}} 1_{\{ X_{2n}=0   \}}  ) 
=\P_o^{\omega}(  Z_n=o) .$$
Donc le lemme \ref{4.5} du point 2, nous donne une borne supérieure du bon 
ordre pour une valeur $\alpha_1$, et le point 1 permet d'étendre cette borne à
tout $\alpha$ de $]0,1[$.
\end{enumerate}
\end{proof}
\section{Borne inférieure de  $\mathbb{E}_0^{\omega}(\alpha^{N_n}) $}
Dans cette partie, on montre la
\begin{prop}  \label{etoh}
 Il existe une constante $c>0$ telle que $Q$  p.s sur   $|\C|=+\infty,$ et pour 
  $n$ assez grand, on ait: 
 $$\E_0^{\omega}(  \alpha^{N_{n}})  \geq  e^{-c n^{\frac{d}{d+2} } }. $$
  \end{prop}

\subsection { Faits généraux pour les marches aléatoires } $\ $ \\
Ce qui suit peut s'appliquer à toute cha\^ine de Markov    $Z$ qui admet
une mesure réversible $m$.  Soit  $A$  un ensemble de points du graphe sur lequel
  $Z$ évolue. On peut écrire
\begin{eqnarray*}
\tilde{\P}^{\omega}_o(Z_{ 2n}=o) &=& \underset{z}{\sum} \tilde{\P}^{\omega}_o(Z_n=z)\tilde{\P}^{\omega}_z
(Z_n=o),\\
&\geq& \underset{z \in A}{\sum}  \tilde{\P}^{\omega}_o(Z_n=z) \times \frac{ m(0)}{m(z)} 
\tilde{\P}^{\omega}_o(Z_n=z),\\
&\geq& \frac{m(0)}{m(A)}[ \underset{z \in A}{\sum} 
\tilde{\P}^{\omega}_o(Z_n=z)]^2,\\
&\geq& \frac{m(0)}{m(A)}[  \tilde{\P}^{\omega}_o(Z_n \in A)]^2.\\
\end{eqnarray*}
Dans notre cas, nous prenons:  $$A=A_r=\{  (a,f) ; \ a\in B_r(\C) \ \text{et
supp}(f)\subset B_r(\C) \}. $$
Rappelons que 
 $m(a,f)= \nu (a) (\frac{1-\alpha}{\alpha})^{ \# \{  i ;
 f(i)=\bar{1} \}  }$.
Donc, on a: 
$$m(A_r) \leq 2d r^d \underset{ k=0...r^d}{\sum} C_{r^d}^{k}
 (\frac{1-\alpha}{\alpha})^k = 
2d \frac{r^d}{ \alpha^{r^d}}.$$ 
Par ailleurs, la structure des ar\^etes sur un produit en couronne, implique
 :
\begin{eqnarray*}
\tilde{\P}^{\omega}_o(Z_n \in A_r) &\geq& \P_0^{\omega}( \forall i \in [|0,n|] \ X_i \in B_r(\C) ) \\
&\geq& \P_0^{\omega}(    \underset{ 0 \leq i \leq n}{\sup} D(0,X_i) \leq r  ),
\end{eqnarray*}
Ainsi, on a:
\begin{eqnarray} \nonumber   \E_0^{\omega}(  \alpha^{N_{2n}}
1_{\{ X_{2n}=0   \}}  ) &  =& 
   \tilde{\P}^{\omega}_o(Z_{ 2n}=o)  \\
   \label{LB} &\geq &
    \frac{ \alpha^{r^d}}{2d r^d} 
[\P_0^{\omega}(    \underset{ 0 \leq i \leq n}{\sup} D(0,X_i) \leq r  ) ]^2 .
\end{eqnarray}
On  est  ramené à trouver une borne inférieure de 
$\P_0^{\omega}(    \underset{ 0 \leq i \leq n}{\sup} D(0,X_i) \leq r  )$, ce qui
est le but de la section suivante.
\subsection{ Borne inférieure de  $\P_0^{\omega} ( \sup_{0 \leq i \leq n } D(0,X_i) \leq r )$} 
Ce type d'estimée inférieure est connue dans le cas de $\Z^d$, mais la 
preuve utilise un principe de  réflexion (voir\cite{SPbook}), que l'on ne peut
utiliser dans  $\C$. Notre démarche  
 (inspirée de  \cite{SP})  pour contourner cette difficulté, 
 utilise 
  des outils analytiques. On prouve:
\begin{prop}\label{borneinfkl} Il existe une constante  $c $ telle que  
$Q$  p.s, pour  $r$ assez grand et pour $n\geq r$, on ait:
$$
 \P_0^{\omega}  (\underset{0 \leq  i \leq  n}{\sup }   D(0,X_i) \leq r  ) \geq
 e^{-c(r +n/r^2)}.$$
\end{prop}
(Pour $ n < r$ cette  probabilité est égale à 1.)
\begin{proof}
Divisons la preuve en 5 temps.
\begin{enumerate} [1.]
\item
Notons $\sigma_r=\min \{j \geq 0 ; \ X_j \not\in B_r(\C) \} $ et considérons
les opérateurs  $P^{B_r(\C)}$ définis sur les fonctions $f:\C\rightarrow \R$
par:
 $$ P^{B_r(\C)} (f)(x) = \ 1_{B_r(\C)} (x) \  
 \E_x^{\omega}(f(X_1) ; \sigma_r >1).$$
  L'opérateur $ P^{B_r(\C)} $ est sous markovien mais
 néamoins symétrique par rapport à la restriction de $\nu$ à
   $B_r(\C)$.
 On note $p^{ B_r ( \C ) } ( , )$ les noyaux de 
 $P^{B_r(\C)} $ et $p^{ B_r ( \C ) }_n( , )$ les noyaux de 
  $(P^{B_r(\C)} )^n$.
   Pour tout  $ x, y \in B_r(\C),$ on a: 
 $$ p^{ B_r ( \C ) } (x,y) =  1_{  \{  (x,y) \in B_r(\C)^2  \}  } 
 \P_x^{\omega} (X_1 = y ) ,$$
 et $$ p^{B_r(\C)}_n(x,y)= 
 \P_x^{\omega}( X_n=y ; \sigma_r   > n) .$$
\item Pour  $n \geq r$ on a successivement:
\begin{eqnarray}
 \P_0^{\omega}(    \underset{ 0 \leq i \leq n}{\sup} D(0,X_i) \leq r  )
\nonumber  &=&  \P_0^{\omega}(  \sigma_r >n)  \\
\nonumber &=& \underset{ y \in B_r(\C) }{\sum}   p^{B_r(\C)}_n(0,y) \\
\nonumber &\geq& \underset{ y \in B_r(\C)  }{\sum}  p^{B_r(\C)}_{r}(0,y) p^{B_r(\C)}_{n-r}(y,y) \\
\nonumber &\geq& e^{-rc}  \underset{ y \in B_r(\C)  }{\sum}   p^{B_r(\C)}_{n-r}(y,y) \\
 \nonumber &\ & \text{(avec }  c=\log(2d ), \ \text{puisque
par connexité,} \\
 \nonumber &\ & \ \text{ il existe un chemin entre } 0 \text{ et } y        ) \\  
 \nonumber &\geq& e^{-rc} Tr(P^{B_r(\C)}_{n-r} )\\
  \label{operat}&\geq& e^{-rc}    (1-\lambda_1)^{n-r},\\ \nonumber
   \end{eqnarray}
o\`u $\lambda_1= \lambda_1 (B_r(\C))$ est la plus petite valeur propre
strictement positive de  l'opérateur  $(Id-P^{B_r(\C)}) 1_{B_r(\C)}.$ 
\begin{eqnarray}
\label{lambda}
\lambda_1 (B_r(\C)) = \underset { \underset{f\neq 0}{\text{supp}(f)
\subset B_r(\C)}}{\inf} \ 
\frac{  \xi(f,f)}{ ||f||_{l^2(\nu)}^2} ,
\end{eqnarray}
avec  
$$\xi(f,f)= \Bigl( (Id-P^{B_r(\C) } ) 1_{B_r(\C)}f|f\Bigr)_{\nu}.$$
Notons que si $\text{supp}(f)\subset B_{r-1}(\C)$, on peut simplifier 
l'expression de  $\xi$.
\begin{eqnarray*}
\xi(f,f)&=& \underset{x\in B_r(\C)}{ \sum} \nu(x)f(x) 
[f(x)- \underset{y\in B_r(\C)}{\sum} p^{B_r(\C)}(x,y)f(y) ]\\
&=& \underset{x\in B_r(\C)}{ \sum} \nu(x)f(x)  
[ \underset{y\in B_r(\C)}{\sum} p^{B_r(\C)}(x,y) (f(x)-f(y))] \\ 
&\ & \hspace{2cm}(\text{Car si } x\in B_{r-1} (\C),\ 
\underset{y\in B_r(\C)}{\sum} p^{B_r(\C)}(x,y)=1. ) \\
&=& \frac{1}{2} \underset{x,y \in B_r(\C) }{\sum} \nu(x)p^{B_r(\C)}(x,y) (f(x)-f(y))^2.
\end{eqnarray*}

\item Prenons maintenant $h(x)= (r-|x| )1_{B_r(\C)}$. On a  
$\text{supp}(h)\subset B_{r-1}(\C), $ donc on peut utiliser la formule
précédente
pour $\xi(h,h)$, et on a:
\begin{enumerate} [$\bullet$]
\item $\xi(h,h) \leq |B_r(\C)|.$
\item $h \geq r/2 $ sur $B_{r/2}(\C)$, donc 
 $||h||_{l^2(\nu)}^2 \geq  (r/2)^2  \nu(B_{r/2}(\C))  $\\
$\hspace*{6.3cm}\geq \frac{1}{2d}(r/2)^2  |B_{r/2}(\C)|.$
\end{enumerate}
Ainsi par (\ref{lambda}), on déduit:
\begin{eqnarray} \label{bornevp}
\lambda_1\leq 8d \frac{ |B_{r} (\C)|}{r^2| B_{r/2} (\C)|}
\end{eqnarray}
\item  L'étape suivante consiste à minorer $|B_{r/2} (\C)|$. On a le lemme
suivant:
\begin{lem} \label{volumekl}
Pour  $p>p_c$, il existe $c>0$ tel que  $Q$ p.s sur  $|\C|=+\infty$, 
 pour  $r$ assez grand, on ait:
$$ |B_r(\C)| \geq  c r^d .$$
\end{lem}
\begin{proof} Remarquons que : $$\text{si} \  m \geq 
\underset{x\in C_r} {\max} D(0,x), \  \text{alors } \ 
\C_r \subset B_m(\C).$$
Maintenant, par  le corollaire 1.3 de  \cite{antal}, il existe une constante
  $\rho'=\rho'(p,d)\geq 1$ telle que   $Q$ p.s, 
$$ \limsup_{N_1(y)\rightarrow +\infty} \frac{D(0,y)}{ N_1(y)} \leq \rho'. $$
Donc, \begin{eqnarray}
\nonumber  
\exists \rho'\geq 1 \ Q \ p.s \ \ \ \exists A_{\omega}\geq 0 
 \ \forall y\in \C ,   \\
 \label{antalbis} \ \ \ \ \ 
( N_1(y)\geq A_{\omega} \Rightarrow D(0,y) \leq 2 \rho N_1(y)).
\end{eqnarray} \\
Soit $r\geq r_{\omega}= \underset{\underset{N_1(z)< A_{\omega}}{z\in \C}}{\max}
D(0,z)$ et prenons   $m=2\rho' r$. Soit maintenant  $x\in \C_r$, on distingue
deux cas:

-soit  $N_1(x)\geq A_{\omega},$ alors par (\ref{antalbis}), $D(0,x)\leq 2 \rho'
N_1(x)\leq 2\rho' r=m$,

-ou bien  $N_1(x) < A_{\omega},$ mais alors par notre choix de  
$r,\ D(0,x)\leq r \leq 2\rho' r=m.$\\

Dans tous les cas, on obtient que pour tout  $r\geq r_{\omega}$ et pour 
tout  $x$ dans $\C_r,\ D(0,x)\leq  m=2\rho n .$ On déduit donc,
 $$\exists \rho' \ \text{tel que } \ Q \ \text{p.s pour $r$ assez grand, } 
 \ \C_r \subset  B_{2\rho' r}(\C).$$
Puis par la propriété   \ref{C^n} (voir Appendice B dans \cite{pierre}),
 on sait qu'il existe une constante $\rho>0$ telle que 
 $Q$ p.s, pour  $r$ assez grand, $|\C_r|\geq \rho r^d.$
Finalement, on a bien l'existence d'une constante  $c>0$ telle  $Q$ p.s  pour
 $r$ assez grand, $ |B_r(\C)| \geq c r^d.$
\end{proof}
\item On déduit alors immédiatement de l'étape 4 et de (\ref{bornevp}) qu'il
 existe  $C>0$  tel que $ Q $ p.s sur  $ |\C|=+\infty, $ 
  pour $ r $ assez grand, 
  $$ \lambda_1(B_r(\C)) \leq  \frac{ C}{ r^2}.$$ 
  Donc pour $r$ assez grand, on peut écrire:
$ \ \ 1-\lambda_1 \geq e^{-2 \lambda_1}$. Et finalement, avec 
  (\ref{operat}), on déduit: 
$$ \ \P_0^{\omega}(    \underset{ 0 \leq i \leq n}{\sup} D(0,X_i) \leq r  ) \geq e^{-c(r + (n-r)/r^2) } \geq 
e^{-c(r + n/r^2) }.$$
\end{enumerate}
\end{proof}
\subsection{ Preuve de la proposition \ref{etoh}.}
Par l'inégalité (\ref{LB}) et par la propriété \ref{borneinfkl}, on déduit, 
qu'il existe $c>0$ tel que $Q$ p.s, pour $r$ assez grand et $n\geq r$, on ait:
\begin{eqnarray*}
\E_0^{\omega}(  \alpha^{N_{2n}}    1_{\{X_{2n}=0}\}) 
 &\geq&     \frac{ \alpha^{r^d}}{2d r^d} e^{-2c(r + n/r^2) } \\
     &\geq&  e^{ -c' ( r^d + n/r^2)    }.
 \end{eqnarray*}
choisissons  $r$  proportionnel à 
 $  n^{\frac{1}{d+2}}$ (qui est bien plus petit que $n$, pour $n$ assez grand),   ainsi on obtient 
 finalement l'existence d'une constante $c>0$ telle que $Q$ p.s sur
 $|\C|=+\infty$, pour $n$ assez grand,
  $$  \E_0^{\omega}(  \alpha^{N_{2n}})  \geq  e^{-c n^{\frac{d}{d+2} } }  . $$
 Pour les temps impairs, on remarque simplement que
  $ \E_0^{\omega}(  \alpha^{N_{2n}}) \leq \E_0(  \alpha^{N_{2n-1}})$. Ce qui 
  achève la preuve de \ref{etoh}.
    
 \section{Questions et  extensions}
 \subsection{Questions ouvertes } $\ $ \\
La première question naturelle, est de savoir si 
 $ \frac{1}{n^{d/d+2}} \log \E_0^{\omega} (\alpha^{N_n}  ) $
 conver-ge lorsque  
 $n$ tend vers l'infini ? Nous savons seulement que cette expression est
 bornée . 
 
 Une autre question intéressante est d'estimer 
  $Q( \E_0^{\omega}(  \alpha^{N_{n}}) )$. Pour la borne supérieure,
  à $n$ fixé on doit estimer la $Q$ probabilité que
  l'inégalité isopérimétrique de la propriété
   \ref{is de base}  ne
  soit pas réalisée. Dans   \cite{barjo}, des calculs assez
  précis donnent un contr\^ole de ce type d'événement  en    $e^{-cn^{\beta}}$ avec  
 $\beta =\frac{d-1}{d+1}<\frac{d}{d+2}$ .
 Or $e^{-cn^{d/d+2}} << e^{-cn^{\beta}}$, et ce terme d'erreur est déjà trop
 "gros" pour la borne espérée.

 \subsection{Extensions} $\ $ \\
La  m\^eme démarche  donne des
bornes supérieures
pour des fonctionnelles plus générales, comme par exemple, 
 $exp^{-\underset{z;L_{z,n}\neq 0}{\sum}F(L_{z,n} ,z)}$ o\`u  $F$ 
 une fonction positive et  $L_{z,n} $  est le nombre de visite de $z$ 
 par la marche $X$ sur l'amas avant le temps  $n$ (voir \cite{mythese}).\\
 \\
 \bf{  Remerciements}\\
 \it{Je remercie  mon directeur de thèse Pierre Mathieu  ainsi que
   Christophe Pittet, Enrique Andjel  et Anna Erschler pour leurs remarques 
 très utiles.}

\nocite{*}
\bibliographystyle{plain}
\bibliography{b}

Clément Rau
\\ \address{Laboratoire d'Analyse, de Topologie et de
Probabilit\'es
\\ Centre de Math\'ematiques et d'Informatique
\\ 39, rue F. Joliot Curie
\\ 13453 Marseille Cedex 13 }
\\  \email{rau@cmi.univ-mrs.fr}
\\ \urladdr{http://www.cmi.univ-mrs.fr/~rau/}

\end{document}